\documentclass[11pt,draft,reqno]{amsart}
\input epsf.tex

\usepackage[all]{xy}
\usepackage{amsthm,array,amssymb,amscd,amsfonts,latexsym}
\xyoption{poly}

\DeclareFontFamily{OT1}{wncyr}{\hyphenchar\font45 }
\DeclareFontShape{OT1}{wncyr}{m}{n}{%
   <5> <6> <7> <8> <9> gen * wncyr
   <10> <10.95> <12> <14.4> <17.28> <20.74>  <24.88>wncyr10}{}
\DeclareFontShape{OT1}{wncyr}{m}{it}{%
   <5> <6> <7> <8> <9> gen * wncyi
   <10> <10.95> <12> <14.4> <17.28> <20.74> <24.88> wncyi10}{}
\DeclareFontShape{OT1}{wncyr}{m}{sc}{%
   <5> <6> <7> <8> <9> <10> <10.95> <12> <14.4>
   <17.28> <20.74> <24.88>wncysc10}{}
\DeclareFontShape{OT1}{wncyr}{b}{n}{%
   <5> <6> <7> <8> <9> gen * wncyb
   <10> <10.95> <12> <14.4> <17.28> <20.74> <24.88>wncyb10}{}
\input cyracc.def
\def\rus{\usefont{OT1}{wncyr}{m}{n}\cyracc\fontsize{9}{11pt}\selectfont}
\def\rusit{\usefont{OT1}{wncyr}{m}{it}\cyracc\fontsize{9}{11pt}\selectfont}
\def\russc{\usefont{OT1}{wncyr}{m}{sc}\cyracc\fontsize{9}{11pt}\selectfont}

\DeclareMathSizes{9}{9}{7}{5} 

\begin{document}
\newtheorem
{theorem}{Theorem}

\newtheorem{proposition}
{Proposition}

\newtheorem{lemma}
{Lemma}
\newtheorem{corollary}
{Corollary}
\newtheorem{question}
{Question}

\theoremstyle{definition}
\newtheorem{definition}
{Definition}
\newtheorem{remark}
{Remark}
\newtheorem{example}
{Example}
\newtheorem{reduction}
{Reduction}
\newtheorem{problem}
{Problem}

\newcommand{\Span}{\operatorname{Span}}
\newcommand{\Symp}{\mbox{\boldmath$\rm Sp$}}
\newcommand{\g}{\mathfrak{g}}
\newcommand{\el}{\mathfrak{l}}
\newcommand{\lt}{\mathfrak{t}}
\newcommand{\lc}{\mathfrak{c}}
\newcommand{\lu}{\mathfrak{u}}
\newcommand{\lr}{\mathfrak{r}}
\newcommand{\Id}{\operatorname{id}}
\newcommand{\id}{\operatorname{id}}
\newcommand{\pr}{\operatorname{pr}}
\newcommand{\Hom}{\operatorname{Hom}}
\newcommand{\sign}{\operatorname{sign}}

\newcommand{\bbA}{{\mathbb A}}
\newcommand{\bbC}{{\mathbb C}}
\newcommand{\bbZ}{{\mathbb Z}}
\newcommand{\bbP}{{\mathbb P}}
\newcommand{\bbQ}{{\mathbb Q}}
\newcommand{\bbG}{{\mathbb G}}
\newcommand{\bbN}{{\mathbb N}}
\newcommand{\bbF}{{\mathbb F}}

\newcommand{\e}{{\lambda}}

\newcommand{\ve}{{\varepsilon}}
\newcommand{\vp}{{\varpi}}

\newcommand{\kbar}{\overline k}

\newcommand{\mo}{\mathopen<}
\newcommand{\mc}{\mathclose>}

\newcommand{\Ker}{\operatorname{Ker}}
\newcommand{\Aut}{\operatorname{Aut}}
\newcommand{\sdp}{\mathbin{{>}\!{\triangleleft}}} 
\newcommand{\Alt}{\operatorname{A}}   
\newcommand{\GL}{\mbox{\boldmath$\rm GL$}}

\newcommand{\PGL}{\mbox{\boldmath$\rm PGL$}}
\newcommand{\SL}{\mbox{\boldmath$\rm SL$}}
\newcommand{\rank}{\operatorname{rank}}
\newcommand{\aut}{\operatorname{Aut}}
\newcommand{\Char}{\operatorname{\rm char\,}} 
\newcommand{\diag}{\operatorname{\rm diag}}
\newcommand{\Gal}{\operatorname{Gal}}
\newcommand{\galois}{\Gal}
\newcommand{\lra}{\longrightarrow}
\newcommand{\SO}{\mbox{\boldmath$\rm SO$}}
\newcommand{\M}{\operatorname{M}}        
\newcommand{\ord}{\mathop{\rm ord}\nolimits}
\newcommand{\Sym}{{\operatorname{S}}}    
\newcommand{\tr}{\operatorname{\rm tr}}
\newcommand{\trace}{\tr}
\newcommand{\ad}{\operatorname{Ad}}

\newcommand{\Res}{\operatorname{Res}}
\newcommand{\Sha}{\mbox{\rus{\fontsize{11}{11pt}\selectfont{SH}}}}
\newcommand{\G}{\mathcal{G}}
\renewcommand{\H}{\mathcal{H}}
\newcommand{\gen}[1]{\langle{#1}\rangle}
\renewcommand{\O}{\mathcal{O}}
\newcommand{\C}{\mathcal{C}}
\newcommand{\Ind}{\operatorname{Ind}}
\newcommand{\End}{\operatorname{End}}
\newcommand{\Spin}{\mbox{\boldmath$\rm Spin$}}
\newcommand{\T}{\mathbf G}
\newcommand{\GT}{\mbox{\boldmath$\rm T$}}
\newcommand{\Inf}{\operatorname{Inf}}
\newcommand{\Tor}{\operatorname{Tor}}
\newcommand{\m}{\mbox{\boldmath$\mu$}}

\newcommand{\A}{{\sf A}}

\newcommand{\D}{{\sf D}}

\newcommand{\Lbd}{{\sf \Lambda}}


\title[Generically multiple transitive
algebraic group actions]{Generically multiple
transitive\\ algebraic group actions}

\author{ Vladimir L. Popov}
\address{Steklov Mathematical Institute,
Russian Academy of Sciences, Gubkina 8, Moscow 119991,
Russia} \email{popovvl@orc.ru}

\thanks{Supported in part by ETH, Z\"urich,
 Switzerland, and Russian grant {\rus
N{SH}--123.2003.01}.}
\date{February 21, 2005}


\begin{abstract}
With every nontrivial connected algebraic group $G$ we
associate a positive integer ${\rm gtd}(G)$ called the
generic transitivity degree of $G$ and equal to the
maximal $n$ such that there is a nontrivial action of
$G$ on an irreducible algebraic variety $X$ for which
the diagonal action of $G$ on $X^n$ admits an open
orbit. We show that ${\rm gtd}(G)\leqslant 2$
(respectively, ${\rm gtd}(G)=\nobreak 1$) for all
solvable (respectively, nilpotent) $G$, and we
calculate ${\rm gtd}(G)$ for all reductive $G$. We
prove that if $G$ is nonabelian reductive, then the
above maximal $n$ is attained for $X=G/P$ where $P$ is
a proper maximal parabolic subgroup of $G$ (but not
only for such homogeneous spaces of $G$). For every
reductive $G$ and its proper maximal parabolic
subgroup $P$, we find the maximal $r$ such that the
diagonal action of $G$ (respectively, a Levi subgroup
$L$ of ~$P$) on $(G/P)^r$ admits an open $G$-orbit
(respectively, $L$-orbit). As an application, we
obtain upper bounds for the multiplicities of trivial
components in some tensor product decompositions. As
another application, we classify all the pairs $(G,
P)$ such that the action of $G$ on $(G/P)^3$ admits an
open orbit, answering a question of {\sc M.~Burger}.
\end{abstract}

\maketitle

 \section {\bf Introduction}
\label{sect.intro}

My starting point was the following question posed to
me by {\sc M.~Burger} in the fall of~2003.

\begin{question} \label{burger} {\it
Let $S$ be a complex connected simple linear algebraic
group and let  $P$ be its proper maximal parabolic
subgroup. For which $S$ and $P$ is there an open
$S$-orbit in $S/P\times S/P\times S/P$?
}
\end{question}

Answering this question led me to the following
general set up. Consider an algebraic action of an
algebraic group $G$ on an algebraic variety $X$. Its
$n$-transitivity means that
\begin{enumerate}
\item[(i)] for the diagonal action of $G$ on
$X^n:=\underbrace{X\times\ldots\times X}_n$, there is
an open $G$-orbit $\mathcal O$; \item[(ii)]
$X^n\setminus\mathcal O$, the boundary of $\mathcal
O$, is the union of ``diagonals" of $X^n$.
\end{enumerate}
All multiple
 transitive (i.e., with $n\geqslant 2$)
actions are classified, \cite{knop}. They constitute a
rather small and not very impressive class: for
$n\geqslant 4$, there are none of them; for $n=3$, it
is only the natural action of ${\bf PGL}_2$ on ${\bf
P}^1$; for $n=2$ and reductive $G$, it is only the
natural action of ${\bf PGL}_{m+1}$ on~${\bf P}^m$.

From the point of view of algebraic transformation
groups, imposing restriction (ii) does not look really
natural. It is more natural to consider the cases
where $G$ acts transitively on all $n$-tuples of
elements of $X$ subject to algebraic inequalities {\it
depending} on the problem under consideration. This
leads to the following definition.

\begin{definition}\label{gentran}
Let $n$ be a positive integer. An algebraic action
$\alpha: G\times X\to X$ of a connected algebraic
group $G$ on an irreducible algebraic variety $X$ is
called {\it generically $n$-tran\-sitive} if the
diagonal action of $\,G$ on $X^n$ is locally
transitive, i.e., admits an open $G$-orbit. If
$\alpha$ is not locally transitive, then $\alpha$ is
called {\it generically $0$-transitive}.
\end{definition}

Below we will see that the class of generically
multiple transitive actions is much more rich and
interesting than that of multiple transitive ones.

Since the projection $X^n\to X^{n-1}, \ (x_1,\ldots,
x_n)\mapsto (x_1,\ldots, x_{n-1})$ is $G$-equivariant,
gene\-ric $n$-transitivity yields generic
$m$-transitivity for $0<m\leqslant n$. As $\dim
X^n=~n\dim X$,
\begin{equation}\label{up}
\text{$\alpha$ is not generically $n$-transitive if
$n\dim X>\dim G$.}
\end{equation}
\begin{definition}\label{defgtdeg}
The {\it generic transitivity degree of an action}
$\alpha$ is
\begin{equation*}
{\rm gtd}(\alpha):={\rm sup}\,n
\end{equation*}
with the supremum taken over all $n$ such that
$\alpha$ is generically $n$-transitive. The {\it
generic transitivity degree of a nontrivial connected
algebraic group} $G$ is
\begin{equation*}
{\rm gtd}(G):= {\rm sup}\;{\rm gtd}\,(\alpha)
\end{equation*}
with the supremum taken over all nontrivial actions
$\alpha$ of $\,G$ on irreducible algebraic varieties.
\end{definition}

Thus ${\rm gtd}(\alpha)$ is a positive integer or
$+\infty$ and the latter holds if and only if $X$ is a
single point. Since every nontrivial group $G$ admits
a nontrivial transitive action, ${\rm gtd}(G)$ is a
well defined positive integer, and by \eqref{up}, we
have
 ${\rm
gtd}(G)\leqslant \dim G$.

Using this terminology, Question~\ref{burger} can be
reformulated as the problem of classifying all proper
maximal parabolic subgroups $P$ of connected simple
algebraic groups $G$ such that the action of $G$ on
$G/P$ is generically $3$-transitive.

\smallskip

In this paper I address the problem of calculating the
generic transitivity degrees of connected linear
algebraic groups. The main results are the following
Theorems~\ref{main1}--\ref{main6}.

\begin{theorem}\label{main1} Let $G$ be
a nontrivial connected linear algebraic group.

\begin{enumerate}
\item[\rm (i)] If $\,G$ is solvable, then $${\rm
gtd}(G)\leqslant 2.$$ \item[\rm (ii)] If $\,G$ is
nilpotent, then $${\rm gtd}(G)=1.$$ \item[\rm (iii)]
If $\,G$ is reductive and $\,\widetilde G\to G$ is an
isogeny, then $${\rm gtd}(\widetilde G)={\rm
gtd}(G).$$ \item[\rm (iv)] If $\,Z$ is a torus and
$S_i$ a connected simple algebraic group, $i=1,\ldots,
d$, then
\begin{equation*}
{\rm gtd}(Z\times S_1\times\ldots\times
S_d)=\underset{i}{\rm max}\,{\rm gtd}(S_i).
\end{equation*}
\item[\rm (v)] If $\,G$ is simple, then ${\rm gtd}(G)$
is given by Table $1$:
\end{enumerate}

\

\vskip 2mm

\begin{center}
\begin{tabular}{c||c|c|c|c|c|c|c|c|c}
\text{ {\rm type of} $G$} & ${\sf A}_l$& ${\sf B}_l,$
$l\geqslant 3$& ${\sf C}_l,$ $l\geqslant 2$& ${\sf
D}_l,$ $l\geqslant 4$& ${\sf E}_6$ & ${\sf E}_7$ &
${\sf E}_8$ &${\sf F}_4$ &${\sf G}_2$
\\[2pt]
 \hline
 &&&&&&&&\\[-9pt]
 ${\rm gtd}(G)$ &$l+2$&$3$&$3$&
$3$ & $4$ &$3$&$2$&$2$&$2$
\\
\end{tabular}
\vskip 4mm
 \centerline{\rm Table $1$}
\end{center}

\end{theorem}

\begin{theorem}\label{main2}
Let $G$ be a connected nonabelian reductive group.
Then there is a proper maximal pa\-ra\-bolic subgroup
$P$ of $G$ such that the generic transitivity degree
of the natural action of $\,G$ on $G/P$ is equal to
${\rm gtd}(G)$.
\end{theorem}

Given two positive integers $i$ and $l$,  $i\leqslant
l$, we put
\begin{gather}\label{Mm}
\begin{gathered}
{\mathcal M}_{li}:=\Bigl\{a\in {\mathbb N} \mid a<
\mbox{\fontsize{14pt}{5mm} \selectfont
$\frac{(l+1)^2}{i(l+1-i)}$}\Bigr\},\\
m_{li}:=\underset{a\in {\mathcal M}_{li}}{\rm max}\,a.
\end{gathered}
\end{gather}

\begin{theorem}\label{main3}
If $G$ is a simple linear algebraic group, then the
generic transitivity degree ${\rm gtd}(G:G/P_{i})$ of
the natural action of $\,G$ on $G/P_i$, where $P_i$ is
a standard proper maximal parabolic subgroup of $\,G$
corresponding to the $i$-th simple root $($see
Section~{\rm \ref{stpar}} below\,$)$, is given by
Table $2$:

\hskip 1mm

\begin{center}
\begin{tabular}{c|cc}
\text{\rm type of} $G$ && ${\rm gtd}(G:G/P_{i})$
\\[2pt]
\hline \hline
&\\[-9pt]
${\sf A}_l$&& $m_{li}$\hskip 2mm {\rm(}\hskip
-1.3mm\mbox{{\fontsize{9pt}{5mm}
 \selectfont\rm see}} \eqref{Mm}\hskip -.5mm{\rm)}
\\[3pt]
 \hline
&\\[-8pt]
 ${\sf B}_l,$ $l\geqslant 3$ &&
 $\begin{matrix}
2
\mbox{\ {\fontsize{9pt}{5mm}
\selectfont\rm if}
$i\neq 1, l,$}\\
3\mbox{\ {\fontsize{9pt}{5mm}
\selectfont\rm if} $i=1, l$}
\end{matrix}$
\\[8pt]
\hline
&\\[-8pt]
${\sf C}_l,$ $l\geqslant 2$&& $\begin{matrix}
2\mbox{\ {\fontsize{9pt}{5mm} \selectfont\rm if}
$i\neq 1, l,$}\\
3\mbox{\ {\fontsize{9pt}{5mm}
\selectfont\rm if} $i=1, l$}
\end{matrix}$
\\[8pt]
\hline
&\\[-8pt]
${\sf D}_l,$ $l\geqslant 4$&& $
\begin{matrix}
2\mbox{\ {\fontsize{9pt}{5mm} \selectfont\rm if}
$i\neq 1, l-1, l,$}\\
3\mbox{\ {\fontsize{9pt}{5mm}
\selectfont\rm if}
$i=1, l-1, l$}
\end{matrix}$
\\[9pt]
\hline
&\\[-8pt]
${\sf E}_6$&&
$\begin{matrix}
\mbox{
\hskip -2mm$4$ {\fontsize{9pt}{5mm} \selectfont\rm if}
$i=1, 6,$}\\
\hskip -1mm 2\mbox{\ {\fontsize{9pt}{5mm}
\selectfont\rm if} $i\neq 1, 6$}
\end{matrix}$
\\[11pt]
\hline
&\\[-8pt]
${\sf E}_7$&&
$\begin{matrix}
2\mbox{\ {\fontsize{9pt}{5mm} \selectfont\rm if}
$i\neq 7,$}\\
3\mbox{\ {\fontsize{9pt}{5mm}
\selectfont\rm if} $i=7$}
\end{matrix}$
\\[8pt]
\hline
&\\[-10pt]
${\sf E}_8$&&
$2$\\[3pt]
\hline
&\\[-8pt]
${\sf F}_4$&&
$2$\\[1pt]
\hline
&\\[-8pt]
${\sf G}_2$&&
$2$\\
\end{tabular}
\vskip 4mm
 \centerline{\rm Table $2$}
\end{center}



\end{theorem}

\vskip 2mm

Since \eqref{Mm} yields $m_{li}\geqslant 3$ for all
$i, l$,  and $m_{li}=3$ if and only if $2i=l+1$, we
obtain

\medskip

\noindent{\bf Corollary~1.} {\it Let $G$ be a
connected simple algebraic group of type ${\sf A}_l$.
Then}
\begin{equation*}
{\rm gtd}(G:G/P_i)\begin{cases}=3
&\mbox{\rm if $2i=l+1$},\\
\geqslant 4&\mbox{\rm otherwise.}\end{cases}
\end{equation*}

\smallskip

The next corollary answers {\sc M.~Bur\-ger}'s
Question~\ref{burger}.

\medskip

\noindent {\bf Corollary~2.} {\it Maintain the
notation of Question $1$. The following properties are
equivalent}:
\begin{enumerate}
\item[\rm (i)] {\it
$(S/P)^3$ contains an open $S$-orbit};
\item[\rm (ii)]
{\it $S$ is of type ${\sf A}_l$, ${\sf B}_l$, ${\sf
C}_l$, ${\sf D}_l$, ${\sf E}_6$ or ${\sf E}_7$, and
$P$ is conjugate to one of the
standard
proper maximal parabolic subgroups $P_i$ of~$\,S$
given by Table $3$:}

 \

\begin{center}
\begin{tabular}{c||c|c|c|c|c|c}
\text{ {\rm type of} $G$} & ${\sf A}_l$& ${\sf B}_l,$
$l\geqslant 3$& ${\sf C}_l,$ $l\geqslant 2$& ${\sf
D}_l,$ $l\geqslant 4$ &${\sf E}_6$ & ${\sf E}_7$
\\[2pt]
 \hline
 &&&&&\\[-9pt]
 $i$ &$1,\ldots, l$&$1, l$&$1,l$&
 $1, l-1, l$& $1, 6$
&$7$
\\
\end{tabular}
\nopagebreak \vskip 4mm
 \centerline{\rm Table $3$}
\end{center}
\end{enumerate}

\begin{remark} Theorem~\ref{main3} shows that, in the
notation of Question~\ref{burger}, there are the cases
where $S$ has an open orbit in $(S/P)^d$ not only for
$d=3$ but also for $d=4$ (for $S$ of types ${\sf A}_l$
with $l\geqslant 2$, and ${\sf E}_6$) and $d\geqslant
5$ (for $S$ of type ${\sf A}_l$ with $l\geqslant 3$).
\end{remark}

\begin{remark}  It would be interesting to
extend Theorem~\ref{main3} by calculating for every
simple $G$ the generic transitivity degree ${\rm
gtd}(G:G/P)$ of the action of $G$ on $G/P$ for every
non-maximal parabolic subgroup $P$. By
Lemma~\ref{ineq} below, if $B$ is a Borel subgroup of
$G$, then ${\rm gtd}(G:G/P)\geqslant {\rm gtd}(G:G/B)$
for every $P$. By Corollary~2 of Proposition~\ref{p-}
below, ${\rm gtd}(G:G/B)=2$ (respectively,~$3$) if $G$
is not (respectively, is) of type ${\sf A}_1$. This
and Theorem~\ref{main3} yield that ${\rm
gtd}(G:G/P)=2$ for every $P$ if $G$ is of types ${\sf
E}_8$, ${\sf F}_4$ or ${\sf G}_2$. \quad $\square$
\end{remark}

According to classical  Richardson's theorem,
\cite{richardson0}, if $P=LU$ is a proper parabolic
subgroup of a connected reductive group $G$ with $U$
the unipotent radical of $P$ and $L$ a Levi subgroup
in $P$, then for the conjugating action there is an
open $P$-orbit in $U$. In general, an open $L$-orbit
in $U$ may not exist. A standard argument (see
Section~\ref{stpar}) reduces the problem of
classifying cases where there is an open $L$-orbit in
$U$ to that for simple $G$. If $P$ is maximal, such a
classification is obtained as a byproduct of our proof
of the above theorems:

\begin{theorem}\label{main4}
Maintain the above notation. Let $P^-$ be a parabolic
subgroup opposite to~$P$.
\begin{enumerate}
\item[\rm (a)] The following properties are
equivalent:
\begin{enumerate}
\item[\rm (i)] $U$ contains an open $L$-orbit,
\item[\rm (ii)] $G/P$ contains an open $L$-orbit,
\item[\rm (iii)] $(G/P)^2\times G/P^-$ contains an
open $G$-orbit.
\end{enumerate}
\item[\rm (b)] If $\,G$ is simple and $P$ is maximal,
then properties {\rm (i)}, {\rm (ii)}, {\rm (iii)}
hold if and only if $P$ is conjugate to one of the
standard maximal parabolic subgroups $P_i$ of $G$
given by Table $4$:

\

\begin{center}
\begin{tabular}{c||c|c|c|c|c|c}
\text{ {\rm type of} $G$} & ${\sf A}_l$& ${\sf B}_l,$
$l\geqslant 3$& ${\sf C}_l,$ $l\geqslant 2$& ${\sf
D}_l,$ $l\geqslant 4$ & ${\sf E}_6$& ${\sf E}_7$
\\[2pt]
 \hline
 &&&&&\\[-9pt]
 $i$ &$1,\ldots, l$&$1, l$&$1,l$&
 $1, l-1, l$&$1, 3, 5, 6$
&$7$
\\
\end{tabular}
\nopagebreak \vskip 4mm
 \centerline{\rm Table $4$}
\end{center}
\end{enumerate}
\end{theorem}

\begin{remark} Finding reductive
subgroups of $G$ that act locally transitively on
$G/P$ is given much attention in the literature, see a
survey in \cite{kime2} and the references therein.
Note that (ii) in Theorem~\ref{main4} is equivalent to
a certain representation theoretic property
 of
the triple $(G, L, P)$ (simplicity of the spectrum),
\cite{vk}. \quad $\square$ \end{remark}

Actually, for maximal $P$, we calculate the generic
transitivity degrees ${\rm gtd}(L:G/P)$ and ${\rm
gtd}(L: U)$ of the $L$-actions on $G/P$ and $U$. In
order to formulate the answer, we put for every two
positive integers $i$ and $l$, $i\leqslant l$,
\begin{gather}\label{Sli}
{\mathcal S}_{li}:=\Bigl\{a\in \{2, 3,\ldots\} \mid
\mbox{\fontsize{14pt}{5mm} \selectfont
$\frac{i}{l+1-i}$}\notin
\Bigl[\mbox{\fontsize{14pt}{5mm} \selectfont
$\frac{a-\sqrt{a^2-4}}{2}$ },
\mbox{\fontsize{14pt}{5mm} \selectfont
$\frac{a+\sqrt{a^2-4}}{2}$ } \Bigr
]\Bigr\}.
\end{gather}
We then have ${\mathcal S}_{li}=\varnothing$ if and
only if $2i=l+1$. For $2i\neq l+1$, we put
\begin{equation}\label{sli}
s_{li}:=\underset{a\in {\mathcal S}_{li}}{\rm max}\,a.
\end{equation}

\begin{theorem}\label{main5}
Maintain the above notation.
\begin{enumerate}
\item[\rm(i)] ${\rm gtd}(L:G/P)={\rm gtd}(L:U)$.
\item[\rm(ii)] If $G$ is simple, then ${\rm
gtd}(L:G/P_i)$ is given by
Table $5$:
\end{enumerate}

\hskip 2mm

\begin{center}
\begin{tabular}{c|c}
\text{\rm type of} $G$ & ${\rm gtd}(L:G/P_{i})$
\\[2pt]
\hline \hline
&\\[-9pt]
${\sf A}_l$& $
\begin{matrix}
1 \mbox{\ {\fontsize{9pt}{5mm} \selectfont\rm if}
$2i=l+1$,}
\\
s_{li} \mbox{\ {\fontsize{9pt}{5mm} \selectfont\rm if}
$2i\neq l+1$\hskip 2mm {\rm(}\hskip
-1.3mm\mbox{{\fontsize{9pt}{5mm}
 \selectfont\rm see}} \eqref{sli},\eqref{Sli}\hskip -.5mm{\rm)}}
\end{matrix}
$
\\[10pt]
 \hline
&\\[-8pt]
 ${\sf B}_l,$ $l\geqslant 3$ &
 $\begin{matrix}
\hskip -1mm 0\mbox{\ {\fontsize{9pt}{5mm}
\selectfont\rm if}
$i\neq 1, l,$}\\
\hskip -2mm 1\mbox{\ {\fontsize{9pt}{5mm}
\selectfont\rm if} $i=1, l$}
\end{matrix}$
\\[8pt]
\hline
&\\[-8pt]
${\sf C}_l,$ $l\geqslant 2$& $\begin{matrix} \hskip
-1mm 0\mbox{\ {\fontsize{9pt}{5mm} \selectfont\rm if}
$i\neq 1, l,$}\\
\hskip -2mm 1\mbox{\ {\fontsize{9pt}{5mm}
\selectfont\rm if} $i=1, l$}
\end{matrix}$
\\[8pt]
\hline
&\\[-8pt]
${\sf D}_l,$ $l\geqslant 4$& $\begin{matrix} 0\mbox{\
{\fontsize{9pt}{5mm} \selectfont\rm if}
$i\neq 1, l-1, l,$}\\
\hskip -1mm 1\mbox{\ {\fontsize{9pt}{5mm}
\selectfont\rm if}
$i=1,$}\\
\hskip 4mm 1\mbox{\ {\fontsize{9pt}{5mm}
\selectfont\rm if} $l$ {\fontsize{9pt}{5mm}
\selectfont\rm
is even and} $i=l-1, l,$}\\
\hskip 1.7mm \mbox{ $2$ {\fontsize{9pt}{5mm}
\selectfont\rm if} $l$ {\fontsize{9pt}{5mm}
\selectfont\rm is odd and} $i=l-1, l$}
\end{matrix}$
\\[23pt]
\hline
&\\[-8pt]
${\sf E}_6$& $\begin{matrix}\hskip -1mm 0\mbox{\
{\fontsize{9pt}{5mm} \selectfont\rm if}
$i=2, 4,$}\\
\hskip -1mm \mbox{$1$
 {\fontsize{9pt}{5mm}
\selectfont\rm if}
$i=3, 5,$}\\
\hskip -2mm \mbox{$2$
 {\fontsize{9pt}{5mm}
\selectfont\rm if} $i=1, 6$}
\end{matrix}$
\\[15pt]
\hline
&\\[-8pt]
${\sf E}_7$& $\begin{matrix}\hskip -1mm 0\mbox{\
{\fontsize{9pt}{5mm} \selectfont\rm if}
$i\neq 7,$}\\
\hskip -2mm 1\mbox{\ {\fontsize{9pt}{5mm}
\selectfont\rm if} $i=7$}
\end{matrix}$
\\[8pt]
\hline
&\\[-10pt]
${\sf E}_8$& $0$\\[3pt]
\hline
&\\[-8pt]
${\sf F}_4$&$0$\\[1pt]
\hline
&\\[-8pt]
${\sf G}_2$&$0$\\
\end{tabular}
\nopagebreak \vskip 4mm
 \centerline{\rm Table $5$}
\end{center}
\end{theorem}

\vskip 2mm

Finally, we show that calculating the generic
transitivity degrees of actions on generalized flag
varieties is closely related to the problem of
decomposing tensor products.

Namely, let $G$ be a connected simply connected
semisimple algebraic group. Fix a Borel subgroup $B$
of $G$ and a maximal torus $T$ of $B$. Let ${\rm
P}_{++}$ be the additive monoid of dominant weights of
$T$ determined by $B$. For $\varpi\in {\rm P}_{++}$,
denote by $E(\varpi)$ a simple $G$-module of highest
weight $\varpi$, and by $\varpi^*$ the highest weight
of the dual $G$-module $E(\varpi)^*$. Let $P(\varpi)$
be the $G$-stabilizer of the unique $B$-stable line in
$E(\varpi)$. The subgroup $P(\varpi)$ of $G$ is
parabolic;
  every parabolic subgroup of $G$ is
obtained this way. If $\varpi$ is fundamental, then
$P(\varpi)$ is maximal.
\begin{theorem}\label{main6} Maintain the above
notation. Let $d$ be a positive integer and let
$\varpi\in {\rm P}_{++}$.
\begin{enumerate}
\item[\rm (i)] if ${\rm gtd}(G:G/P(\varpi))\geqslant
d$,
then
\begin{equation}\label{<2}
\dim \bigl(E(n_1\varpi^*)\otimes\ldots \otimes
(E(n_d\varpi^*)\bigr)^G\leqslant 1\quad\mbox{for every
$n_1,\ldots, n_d\in \mathbb Z_+$};
\end{equation}
\item[\rm(ii)] if $\varpi$ is fundamental, the
converse it true, i.e., {\rm\eqref{<2}} yields ${\rm
gtd}(G:G/P(\varpi))\geqslant d$.
\end{enumerate}
\end{theorem}

As an application of Theorems~\ref{main6},
\ref{main3}, we obtain upper bounds of
the multiplicities of trivial components in some
tensor product decompositions.

\begin{example} Let $\varpi_i$ be the $i$th fundamental
weight of $G={\bf SL}_m$. Then
\begin{equation*}
\dim \bigl(E(n_1\varpi_i)\otimes\ldots \otimes
(E(n_d\varpi_i)\bigr)^G\leqslant 1\quad\mbox{for every
$n_1,\ldots, n_d\in \mathbb Z_+$}
\end{equation*}
if and only if $d<m^2/(im-i^2)$.
\end{example}

In \cite{popov} we develop further the latter topic.

\vskip 2.5mm

\noindent{\it Acknowledgments.} I am grateful to {\sc
M.~Burger} for posing Question~\ref{burger}. My thanks
go to {\sc V.~Kac},  {\sc D.~Shmel'kin}, and {\sc
J.~Weyman} for useful information, and to {\sc
Z.~Reichstein} for remarks. I am especially indebted
to {\sc A.~Schofield} who communicated me the argument
used in the proof of Theorem~\ref{main3}, and to the
referee for pointing out some gaps in the first
version of the paper.

\section{\bf Conventions, notation, and terminology}

In this paper the characteristic of the base
algebraically closed field $k$ is equal to $0$.  The
reason is that in some of the proofs I use the
classification results from \cite{kimu1}, \cite{kkiy},
\cite{kkti} that are obtained under this constraint on
${\rm char}\,k$. The problem of extending the main
results to positive characteristic (perhaps with some
small primes excluded) looks manageable. I did not
attempt to solve it here.

Below every action of algebraic group is algebraic
(morphic).

Given actions of $\,G$ on $X_1,\ldots, X_n$, the
action of $\,G$ on $X_1\times\ldots\times X_n$ means
the diagonal action.

Given the subgroups $S$ and $H$ of $G$, the action of
$S$ on $G/H$ is always the natural action induced by
left translations. It is denoted by $(S : G/H)$.

Let $P=LU$ be a parabolic subgroup of $G$ with $U$ the
unipotent radical and $L$ a Levi subgroup in $P$, and
let $\mathfrak u$ be the Lie algebra of $U$. Then $(L
: U)$ and $(L : \mathfrak u)$ denote  the conjugating
and adjoint actions of $L$ on $U$ and $\mathfrak u$
respectively.

Given an action of a group $G$ on a set $X$, we denote
by $G\cdot x$ and $G_x$ respectively the $G$-orbit and
$G$-stabilizer of a point $x\in X$. The fixed point
set $G$ on $X$ is denoted by $X^G$.

Throughout the paper the notion of {\it general point}
is used. By that it is meant a point lying off a
suitable closed subvariety.

$G^0$ is the identity component of an algebraic group
$G$.

$(G, G)$ is the commutator group of a group $G$.

Given a root system with a base
$\Delta=\{\alpha_1,\ldots,\alpha_r\}$, we enumerate
the simple roots $\alpha_1,\ldots,\alpha_r$ as in
\cite{bourbaki}. If $\varpi_1,\ldots,\varpi_r$ is the
system of fundamental weights corresponding to
$\Delta$ and $\lambda=a_1\varpi_1+\ldots+a_r \varpi_r$
is a weight, then we write $\lambda=(a_1,\ldots,
a_r)$. The labelled Dynkin diagram of $\Delta$, where
the label of $i$-th vertex is $a_i$, will be called
the {\it Dynkin diagram of} $\lambda$ (if $a_i=0$,
then the $i$-th label is dropped). Note that if
$\lambda=c_1\alpha_1+\ldots+ c_r\alpha_r$, then
$(a_1,\ldots,a_r)$ is the linear combination with
coefficients $c_1,\ldots, c_r$ of the rows of the
Cartan matrix of $\Delta$.

We call a connected linear algebraic group {\it
simple} if it has no proper closed normal subgroups of
positive dimension.

$k^n$ is the $n$-dimensional coordinate space of
column vectors over $k$.

$k[X]$ is the algebra of regular functions of an
algebraic variety $X$.

$k(X)$ is the field of rational functions of an
irreducible algebraic variety $X$.

We put $\mathbb Z_+=\{0, 1, 2,\ldots\}$, $\mathbb
N=\{1, 2,\ldots \}$, and $[a, b]=\{x\in \mathbb R\mid
a\leqslant x \leqslant b\}$.

\section{\bf Properties of
generically multiple transitive actions}

In this section we establish some general properties
of generically multiple transitive actions that will
be used in the proof of
Theorems~\ref{main1}--\ref{main4}.

Let $\sigma: G\to H$ be a surjective homomorphism of
nontrivial connected algebraic groups.  Let $\alpha $
be an action of $H$ on an irreducible variety $X$, and
let ${}^\sigma\!\alpha$ be the action of $G$ on $X$
defined by
\begin{equation}\label{br}
{}^\sigma\!\alpha(g, x):=\alpha(\sigma(g), x).
\end{equation}

\begin{lemma}\label{>}
Maintain the above notation. Then
\begin{enumerate}
\item[\rm(i)] $ {\rm gtd}({}^\sigma\!\alpha)={\rm
gtd}(\alpha), $ \item[\rm(ii)] ${\rm gtd} (G)\geqslant
{\rm gtd} (H).$
\end{enumerate}
\end{lemma}
\begin{proof} (i) is clear, and (ii) follows from (i)
and Definition~\ref{defgtdeg}. \quad $\square$
\renewcommand{\qed}{}
\end{proof}

\begin{lemma} \label{ineq}
Let $\alpha_i$ be an action of a connected algebraic
group $G$ on an irreducible algebraic variety $X_i$,
$i=1, 2$. Assume that there exists a $G$-equivariant
dominant rational map $\varphi: X_1\dasharrow X_2$.
\begin{enumerate}
\item[\rm (i)] ${\rm gtd}(\alpha_1)\leqslant {\rm
gtd}(\alpha_2)$. \item[\rm(ii)] If $\,\dim X_1=\dim
X_2$, then ${\rm gtd}(\alpha_1)={\rm gtd}(\alpha_2)$.
\end{enumerate}
\end{lemma}

\begin{proof} (i) Assume that the
action of $\,G$ on $X_1^n$ is locally transitive.
Since the rational map $\varphi^n: X_1^n\dasharrow
X_2^n, \ (x_1,\ldots, x_n)\mapsto
(\varphi(x_1),\ldots, \varphi(x_n)),$ is
$G$-equivariant and dominant, the indeterminacy locus
of $\varphi^n$ lies in the complement to the open
$G$-orbit in $X_1^n$, and the image of this orbit
under $\varphi^n$ is a $G$-orbit open in $X_2^n$.
Definitions~\ref{gentran},~\ref{defgtdeg} now yield
that (i) holds.

(ii) Assume that $\,\dim X_1=\dim X_2$ and there is an
open $\,G$-orbit $\mathcal O$ in $X_2^n$. Since
$\varphi^n$ is $G$-equivariant and dominant, there is
a point $z\in X_1^n$ such that the orbit $G\cdot z$
lies off the indeterminacy locus of $\varphi^n$ and
$\varphi^n(G\cdot z)=\mathcal O$. This yields  $\dim
X_2^n=\dim X_1^n\geqslant \dim G\cdot z\geqslant \dim
\mathcal O=\dim X_2^n$. Hence $\dim X_1^n=\dim G\cdot
z$, i.e., $G\cdot z$ is open in $X_1^n$. From (i) and
Definitions~\ref{gentran},~\ref{defgtdeg} we now
deduce that (ii) holds. \quad $\square$
\renewcommand{\qed}{}
\end{proof}

\noindent{\bf Corollary.}\label{biratinv} {\it The
generic transitivity degree is a birational invariant
of actions, i.e., if $\varphi$ in Lemma {\rm
\ref{ineq}} is a birational isomorphism, then ${\rm
gtd}(\alpha_1)= {\rm gtd}(\alpha_2)$.}\quad $\square$
\renewcommand{\qed}{}

\begin{lemma}\label{actproduct}
Let a connected algebraic group $G_i$ act on an
irreducible variety $X_i$, $i=1,\ldots, d$. All these
actions are generically $n$-transitive if and only if
the natural action of $\,G_1\times\ldots\times G_d$ on
$X_1\times\ldots\times X_d$ is generically
$n$-transitive.
\end{lemma}

\begin{proof} This easily follows from
Definition~\ref{gentran}. \quad $\square$
\renewcommand{\qed}{}
\end{proof}

\begin{lemma} \label{reduction} Let
an algebraic group $G$ act on irreducible algebraic
varieties $X$ and $Y$.
\begin{enumerate}
\item[\rm (a)] The following properties are
equivalent:
\begin{enumerate}
\item[\rm (i)] The action of $\,G$ on $X\times Y$ is
locally transitive. \item[\rm (ii)]  The action of
$\,G$ on $X$ is locally transitive and if $\,H$ is the
$G$-stabilizer of a general point of $X$, then the
natural action of $\,H$ on $Y$ is locally transitive.
\end{enumerate}
\item[\rm (b)] Assume that {\rm (i)} and {\rm (ii)}
hold. If $x\in X$ and $y\in Y$ are the points such
that the orbits $G\cdot x$ and $G_x\cdot y$ are open
 in $X$ and $Y$ respectively, then the orbit $G\cdot
(x, y)$ is open in $X\times Y$.
\end{enumerate}
\end{lemma}
\begin{proof} Note that for every
point $z=(x, y)\in X\times Y$, we have
\begin{equation}\label{stab}
G_z=(G_x)_z=(G_x)_y.
\end{equation}
Assume that (i) holds and the orbit $G\cdot z$ is open
in $X\times Y$. As the natural projection $\pi^{}_X:
X\times Y\to X$ is $G$-equivariant and surjective, the
orbit $G\cdot x=\pi^{}_X(G\cdot z)$ is open in~$X$.
The fiber $\pi_X^{-1}(x)$ is $G_x$-stable, contains
$z$, and $G_x\cdot z=G\cdot z\cap \pi_X^{-1}(x)$.
Hence the orbit $G_x\cdot z$ is open in
$\pi_X^{-1}(x)$. As the natural projection
$\pi_X^{-1}(x)\to Y$ is a $G_x$-equivariant
isomorphism, (ii) and (b)  follow from \eqref{stab}.

Conversely, assume that (ii) holds. Then, if the above
point $z$ is general, \eqref{stab} yields $\dim G\cdot
z=\dim G-\dim G_z=\dim G-\dim (G_x)_y= (\dim G-\dim
G_x)+(\dim G_x-\dim (G_x)_y)=\dim X+\dim Y=\dim
(X\times Y)$, i.e., $G\cdot z$ is open in $X\times Y$.
So (i) holds.
 \quad $\square$
\renewcommand{\qed}{}
\end{proof}

\begin{lemma}\label{parabolic}
 Let $P$ be a proper parabolic subgroup of a
connected reductive group $G$. Let $w_0$ be the
element of the Weyl group of $\,G$ with maximal length
$($with respect to a system of simple reflections$)$,
and let $\overset{.}w_0$ be a representative of $w_0$
in the normalizer of a maximal torus of $P$. Let $x\in
G/P$ be a point corresponding to the coset $P$. Then
the $G$-orbit of point $z:=(x, \overset{.}w_0\cdot x)$
is open in $(G/P)^2$.
\end{lemma}
\begin{proof} It follows from Bruhat decomposition,
cf.\;\cite[14.12, 14.14]{borel}, that
$P\overset{.}w_0P$ is open in $G$. Hence the $P$-orbit
of $\overset{.}w_0\cdot x$ is open in $G/P$. Since
$G_x=P$, the claims now follows from Lemma
\ref{reduction} and \eqref{stab}. \quad $\square$
\renewcommand{\qed}{}
\end{proof}
\noindent{\bf Corollary.} {\it Let $\alpha$ be the
action of a connected linear algebraic group $G$ on
$G/P$, where $P$ is a proper parabolic subgroup of
$\,G$. Then ${\rm gtd}(\alpha)\geqslant 2$.}

\begin{proof} Since $P$ contains the
radical of $\,G$, we may assume that $G$ is reductive
and then apply Lemma~\ref{parabolic}.\quad $\square$
\renewcommand{\qed}{}
\end{proof}

\begin{proposition}\label{>1}
Let $G$ be a nontrivial connected linear algebraic
group.
\begin{enumerate}
\item[(i)] If $\,G$ is nonsolvable, then ${\rm
gtd}(G)\geqslant 2$. \item[(ii)] If $\,G$ is solvable,
then ${\rm gtd}(G)\leqslant 2$. \item[(iii)] If $\,G$
is nilpotent, then ${\rm gtd}(G)=1$.
\end{enumerate}
\end{proposition}
\begin{proof}
Since $G$ is nonsolvable if and only if it contains a
proper parabolic subgroup, \cite[6.2.5]{springer},
statement (i) follows from the Corollary of
Lemma~\ref{parabolic} and Definition~\ref{defgtdeg}.

Assume now that $G$ is solvable (respectively,
nilpotent). Let $\alpha$ be an action of $G$ on an
irreducible algebraic variety $X$ such that
\begin{equation}\label{assumpt}
{\rm gtd}(G)={\rm
gtd}(\alpha).
\end{equation}
Since ${\rm gtd}(G)\geqslant 1$, there is an open
$G$-orbit in $X$. By the Corollary of
Lemma~\ref{ineq}, we may, maintaining equality
\eqref{assumpt}, replace $X$ by this orbit. So we may
(and shall) assume that $\alpha$ is the action of
$\,G$ on $X=G/H$ where $H$ is a proper closed subgroup
of~$\,G$. If $Q$ is a proper maximal closed subgroup
of $\,G$ containing $H$, then the existence of the
natural $G$-equivariant morphism $G/H\to G/Q$,
Lemma~\ref{ineq}(i), Definition~\ref{defgtdeg}, and
equality \eqref{assumpt} yield that we may,
maintaining \eqref{assumpt}, replace $H$ by $Q$. So we
may (and shall) assume that $H$ is a proper maximal
closed subgroup of $\,G$. If $H$ contains a nontrivial
closed normal subgroup $N$ of $G$, then $N$ acts
trivially on $X$, so maintaining
 the generic transitivity degree
of the action on $X$ and the assumption that $G$ is
solvable (respectively, nilpotent), the group $G$ may
be replaced by $G/N$. So we may (and shall) also
assume that $H$ contains no nontrivial closed normal
subgroups of~$G$. Finally, by Lemma~\ref{ineq}(ii),
replacing $H$ by $H^0$, we may (and shall) assume that
$H$ is connected.

Use now that in every nontrivial connected nilpotent
linear algebraic group the dimension of center is
positive, and the dimension of normalizer of every
proper subgroup is strictly bigger than that of this
subgroup, cf.\;\cite[17.4]{humphreys}. From this and
the properties of $H$ we deduce that if $G$ is
nilpotent, then $\dim G=1$ and $H$ is trivial.
By~\eqref{up}, we then obtain ${\rm gtd}(G)=1$. This
proves (iii).

Assume now that $G$ is not nilpotent (but solvable).
Then the set $G_u$ of all unipotent elements of $G$ is
a nontrivial closed connected normal subgroup of $G$,
and if $T$ is a maximal  torus of $G$, then $G$ is a
semidirect but not direct product of $T$ and $G_u$,
see \cite[6.3]{springer}. Since $H$ is a connected
solvable group as well, if $S$ is maximal torus of
$H$, then $H$ is a semidirect product of $S$ and
$H_u$. We have $H_u\subseteq G_u$ and we may (and
shall) take $S$ and $T$ so that $S\subseteq T$. Since
$H$ contains no nontrivial closed normal subgroups of
$G$, we have $H_u\neq G_u$. Therefore if $S\neq T$,
then $SG_u$ is a proper closed subgroup of $G$
containing $H$. As $H$ is maximal, this is impossible.
Thus~$S=T$.

Consider now the center $Z_{G_u}$ of $G_u$. By
\cite[6.3.4]{springer}, we have $\dim Z_{G_u}\geqslant
1$. Since $Z_{G_u}$ is a closed normal subgroup of
$G$, it is stable with respect to the conjugating
action of $T$. Since $Z_{G_u}$ is a commutative
unipotent group, it is $T$-equivariantly isomorphic to
the vector group of its Lie algebra $\mathfrak
z_{G_u}$ on which $T$ acts via the adjoint
representation (if we embed $G$ in some ${\bf GL}_n$,
the exponential map is a $T$-equivariant isomorphism
of the last group to $Z_{G_u}$). As the $T$-module
$\mathfrak z_{G_u}$ is a direct sum of one-dimensional
submodules, from this we deduce that $Z_{G_u}$
contains a one-dimensional  closed subgroup $U$
normalized by $T$.

Since $H=TH_u$, the subgroup $U$ is normalized by $H$
as well. Note now that the subgroup $H_u\cap Z_{G_u}$
is trivial. Indeed, since $T$ normalizes $H_u$ and
$G=TG_u$, this subgroup is normal in $G$, and as $H$
contains no nontrivial closed normal subgroups of $G$,
the claim follows. From this we deduce that the
subgroup $H_u\cap U$ is trivial. Since $\dim U=1$ and
$U$ is normalized by $H$, this easily yields that $HU$
is a closed subgroup of dimension $\dim H +1$. The
maximality of $H$ then yields $HU=G$. We now conclude
that the variety $X=G/H$ is isomorphic to the affine
line~${\bf A}^1$.

It is well known that ${\rm Aut}\,{\bf A}^1$ coincides
with the group ${\bf Aff}_1$ of all affine
transformations of~${\bf A}^1$. Hence the action of
$G$ on $X$ induces a homomorphism $\varphi: G\to {\bf
Aff}_1$. We have ${\rm ker} \varphi \subseteq H$.
Since $H$ contains no nontrivial closed normal
subgroups of $G$, this yields that $\varphi$ is
injective. As ${\bf Aff}_1$ is a connected
$2$-dimensional linear algebraic group, every its
proper subgroup is abelian, hence nilpotent. Since $G$
is not nilpotent, from this we deduce that $\varphi$
is an isomorphism. It now remains to note that using
Lemma~\ref{reduction}(a) and \eqref{up} one easily
verifies that the the generic transitivity degree of
the action of ${\bf Aff}_1$ on ${\bf A}^1$ is equal to
$2$. This proves (ii). \quad $\square$
\renewcommand{\qed}{}
\end{proof}

Let $P=LU$ be a proper parabolic subgroup of a
connected reductive group $G$ with $U$ the unipotent
radical of $P$ and $L$ a Levi subgroup in $P$. Let
$P^-\hskip -1mm$ be the unique
(cf.~\cite[14.21]{borel}) parabolic subgroup opposite
to $P$ and containing $L$, and let $U^-$ be the
unipotent radical of $P^-$. Denote by $p\in G/P$ and
$p^-\in G/P^-$ the points corresponding
 to the cosets $P$ and $P^-$. The orbits
$U^-\cdot p$ and $U\cdot p^-$ are $L$-stable and open
in $G/P$ and $G/P^-$, \cite[14.21]{borel}. Let
$\mathfrak u$ and $\mathfrak u^-$ be the Lie algebras
of $U$ and $U^-$.

\begin{proposition}\label{p-} Maintain the above
notation.
\begin{enumerate}
\item[\rm(i)]  $U$, $\mathfrak u$, $U\cdot p^-$  are
$L$-isomorphic varieties; the same holds for $U^-$,
$\mathfrak u^-$, $U^-\cdot p$. 
So, ${\rm gtd}(L:U)={\rm gtd}(L:\mathfrak u),\ {\rm
gtd}(L:U^-)={\rm gtd}(L:\mathfrak u^-). $
\item[\rm(ii)] $(L:\mathfrak u^-) =
{}^\sigma\!(L:\mathfrak u)$ {\rm(}see \eqref{br}{\rm)}
for some  $\sigma\in {\rm Aut}\,L$, and ${\rm gtd}(L :
\mathfrak u)={\rm gtd}(L : \mathfrak u^-). $
\item[\rm(iii)] ${\rm gtd}(L:U)={\rm gtd}(L:G/P)={\rm
gtd}(L:U^-)={\rm gtd}(L:G/P^-).$
 \item[\rm(iv)]
$(G/P)^d\times G/P^-$ contains an open $G$-orbit if
and only if $d\leqslant 1+{\rm gtd}(L:U^-).$
\end{enumerate}
\end{proposition}
\begin{proof} (i) We may assume that $G\subseteq {\bf GL}_n$,
\cite[2.3.7]{springer}. Then elements of $U$ and
$\mathfrak  u$ are respectively unipotent and
nilpotent matrices. Hence the $L$-equivariant
morphisms $\mathfrak u\to U$, \ $Y\mapsto
\sum_{i=0}^{n-1}\frac{1}{i!}Y^i$, and
 $U\to \mathfrak  u$,\
$Z\mapsto -\sum_{i=1}^{n-1}\frac{1}{i}(I_n-Z)^i$ are
inverse to one another. Thus $U$ and $\mathfrak u$ are
$L$-isomorphic. Since $U\rightarrow U\cdot p^-$,
$u\mapsto u\cdot p^-$, is an $L$-equivariant
isomorphism, cf.\;\cite[14.21]{borel}, $U$ is
$L$-isomorphic to $U\cdot p^-$. For $U^-$, $\mathfrak
u^-$, $U^-\cdot p$ the arguments are the same.

(ii) The first claim follows from the fact that $L$ is
a connected reductive group and the $L$-modules
 $\mathfrak u$ and
$\mathfrak u^-$ are dual to one another (cf., e.g.,
\cite{roh} and Section~\ref{stpar}). The second
 follows from the first
 and Lemma \ref{>}.

(iii) Since $U^-\cdot p$ is open in $G/P$, we deduce
from the Corollary of Lemma~\ref{ineq} and (i) that
${\rm gtd}(L:G/P)= {\rm gtd}(L:U^-)$. Analogously,
${\rm gtd}(L:G/P^-)= {\rm gtd}(L:U)$. The claim now
follows from (i) and (ii).

(iv)  Since $G_{p^-}=P^-$ and $P^-\cdot p=U^-L\cdot
p=U^-\cdot p$, the orbit $G_{p^-}\cdot p$ is open in
$G/P$.  Lemma~\ref{reduction} then yields that  the
orbit $G\cdot z$, where $z=(p, p^-)\in G/P\times
G/P^-$, is open in $G/P\times G/P^-$. Since $G_z=P\cap
P^-=L$, applying Lemma~\ref{reduction} again, we
deduce that $(G/P)^d\times G/P^-$ for $d\geqslant 2$
contains an open $G$-orbit if and only if
$(G/P)^{d-1}$ contains an open $L$-orbit. By
Definition~\ref{defgtdeg}, the latter holds if and
only if $d-1\leqslant {\rm gtd}(L:G/P)$. On the other
hand, since $U^-\cdot p$ is open in $G/P$, the
Corollary of Lemma~\ref{ineq} and (i) yield that ${\rm
gtd}(L:G/P)= {\rm gtd}(L:U^-)$.
  \quad
$\square$
\renewcommand{\qed}{}
\end{proof}
\begin{corollary}\label{cor}
\begin{enumerate}
\item[\rm (i)] ${\rm gtd}(G:G/P)\geqslant 1+{\rm
gtd}(L:\mathfrak u^-)$. \item[\rm (ii)] If $P$ and
$P^-$ are conjugate, then ${\rm gtd}(G:G/P)= 2+{\rm
gtd}(L:\mathfrak u^-)$.
\end{enumerate}
\end{corollary}
\begin{proof} This follow
from Definition~\ref{defgtdeg}, Proposition~\ref{p-},
and two remarks: (a) As the projection $(G/P)^d\times
G/P^-\rightarrow (G/P)^d$ is $G$-equivariant, the
existence of an open $G$-orbit in $(G/P)^d\times
G/P^-$ yields its existence in $(G/P)^d$. (b) The
assumption in (ii) yields $(G/P)^{d+1}=(G/P)^d\times
G/P^-$.\quad $\square$
\renewcommand{\qed}{}
\end{proof}
\begin{corollary}\label{abelianrad} Let $P$
be conjugate to $P^-$. The following properties are
equivalent:
\begin{enumerate}
\item[\rm(i)] $(G:G/P)$ is generically $3$-transitive,
  \item[\rm (ii)]
$(L:\mathfrak u^-)$ is locally transitive.\quad
$\square$
\renewcommand{\qed}{}
\end{enumerate}
\end{corollary}

\begin{remark} \label{abel}
If $U$ is abelian, then, by \cite{richardson}, the
number of $L$-orbits in $U$ is finite, hence, by
Proposition 2, $(L:\mathfrak u^-)$ is locally
transitive (in this case $U^-$ is abelian as well, see
Section~\ref{stpar}, and, by Proposition 2, the number
of  $L$-orbits in $\mathfrak u^-$ is equal to that in
$\mathfrak u$). The parabolic subgroups $P$ whose
unipotent radical is abelian are easy to classify,
\cite{rrs}. Every such $P$ is maximal. If $G$ is
simple, then up to conjugacy all such $P$ are
exhausted by the following standard parabolic
subgroups $P_i$:

\vskip 3mm

\begin{center}
\begin{tabular}{c||c|c|c|c|c|c}
\text{type of $G$} & ${\sf A}_l$& ${\sf B}_l$& ${\sf
C}_l$& ${\sf D}_l$& ${\sf E}_6$ & ${\sf E}_7$
\\[2pt]
 \hline
 &&&&&&\\[-9pt]
 $i$ &$1,\ldots,
l$&$1$&$l$&$1, l-1, l$&$1, 6$&$7$
\\
\end{tabular}
\end{center}
\end{remark}

\begin{corollary} Let $G$ be
a connected semisimple group and let $B$ be a Borel
subgroup of $G$.
\begin{enumerate}
\item[\rm(a)] $2\leqslant{\rm gtd}(G:G/B)\leqslant 3$.
\item[\rm(b)] The following properties are equivalent:
\begin{enumerate}
\item[\rm(i)] ${\rm gtd}(G:G/B)=3$,   \item[\rm (ii)]
$G$ is locally isomorphic to ${\bf
SL}_2\times\ldots\times{\bf SL}_2$.
\end{enumerate}
\end{enumerate}
\end{corollary}

\begin{proof} Use the above notation for $P=B$.
By Lemma~\ref{parabolic}, we have ${\rm
gtd}(G:G/B)\geqslant 2$.

Assume that $(G:G/B)$ is generically 3-transitive.
Since $P$ is conjugate to $P^-$, and $L$ is a maximal
torus of $G$, Corollary~\ref{abelianrad} yields $\dim
L\geqslant \dim U^-=$ num\-ber of negative roots
$\geqslant$  number of simple roots $=\dim L$.
Therefore every positive root is simple, whence $G$ is
locally isomorphic to ${\bf
SL}_2\times\ldots\times{\bf SL}_2$. Now
Proposition~\ref{isogen} below and
Lemma~\ref{actproduct} yield that replacing $G$ by
${\bf SL}_2$ does not change $(G:G/B)$. For $G={\bf
SL}_2$, we have $G/B={\bf P}^1$. It is classically
known that the natural action of ${\bf PGL}_2={\rm
Aut}\,{\bf P}^1$ on ${\bf P}^1$ is $3$-transitive but
not $4$-transitive. Whence ${\rm gtd}(G:G/B)=3$. \quad
$\square$
\renewcommand{\qed}{}
\end{proof}

\begin{lemma}\label{excl-reduct} Let $\alpha$
be the action of a connected reductive group $G$ on
$G/H$, where $H$ is a proper closed reductive subgroup
of $\,G$. Then ${\rm gtd}(\alpha)=1$.
\end{lemma}
\begin{proof} Assume the contrary. Then,
since $H$ is the $G$-stabilizer
of a point of $\,G/H$,
 Lemma~\ref{reduction} yields that $G/H$ contains
  an open
$H$-orbit. On the other hand, by \cite{luna},
\cite{nisn}, reductivity of $H$ and $G$ yields that
the action of $H$ on $G/H$ is stable, i.e., a general
$H$-orbit is closed in $G/H$. Hence the action of $H$
on $G/H$ is transitive. But this is impossible since
the fixed point set of this action is nonempty.\quad
$\square$
\renewcommand{\qed}{} \end{proof}

\begin{proposition}\label{1-2} Let $G$ be
a nontrivial connected reductive group. Then
\begin{enumerate}
\item[\rm (i)] ${\rm gtd}(G)=1$ if and only if $\,G$
is abelian $($i.e., a torus$)$; \item[\rm (ii)] if
$\,G$ is nonabelian, then for some proper maximal
parabolic subgroup $P$ of $\,G$,
\begin{equation}\label{alphaG}
{\rm gtd}(G)={\rm gtd}(\alpha),
\end{equation}
where $\alpha$ is the action of $\,G$ on $G/P$.
\end{enumerate}
\end{proposition}

\begin{proof} Since solvable reductive groups are
tori, (i) follows from Proposition~\ref{>1}.

Assume now that $G$ is nonabelian. The same argument
as in the proof of Proposition~\ref{>1} shows that
there a proper closed maximal subgroup $H$ of $G$ such
that for the action $\alpha$ of $G$ on $G/H$ condition
\eqref{alphaG} holds.
 It is known, cf.\;\cite[30.4]{humphreys},
that the maximality of $H$ yields that $H$  is either
a reductive or a parabolic subgroup of $\,G$. But (i)
yields ${\rm gtd}(G)\geqslant 2$. Therefore
\eqref{alphaG} and Lemma~\ref{excl-reduct} rule out
the first possibility. Thus $H$ is a proper maximal
parabolic subgroup of $\,G$. This proves (ii). \quad
$\square$
\renewcommand{\qed}{}
\end{proof}

\begin{remark}
It would be interesting to classify subgroups $Q$ of
connected nonabelian reductive groups $G$ such that
the action of $\,G$ on $G/Q$ is generically
$2$-transitive. By Lemma~\ref{reduction}, this
property is equivalent to the existence of an open
$Q$-orbit in $G/Q$. Lemma~\ref{excl-reduct} shows that
such $Q$ is not reductive. By Lemma~\ref{parabolic},
every proper parabolic subgroup of $\,G$ has this
property. However the following example shows that
nonparabolic subgroups with this property exist as
well.

\begin{example}\label{nonred} Take a  reductive group $G$ such
that the longest element $w_0$ of the Weyl group of
$\,G$ (with respect to a system of simple reflections)
is not equal to $-{\rm id}$. Let $B$ be a Borel
subgroup of $\,G$ with the unipotent radical $U$. Fix
a maximal torus $T$ of $B$ and let $\overset{.}{w}_0$
be a representative of $w_0$ in the normalizer of $T$.
The above condition on $w_0$ yields that $T$ contains
a subtorus $T'$ of codimension $1$ that is not stable
with respect to the conjugation by $\overset{.}{w}_0$.
Then $T=\overset{.}{w}^{}_0T'\overset{.}{w}_0^{-1}T'$
because of the dimension reason. Since
$\overset{.}{w}^{}_0B\overset{.}{w}_0^{-1}B$ is open
in $G$, this yields that $Q:=T'U$ is a nonparabolic
subgroup of $\,G$ such that the action of $\,G$ on
$G/Q$ is $2$-transitive. The generalizations of this
construction replacing $1$ by a bigger codimension and
$B$ by a parabolic subgroup are clear.
\end{example}
\end{remark}

\begin{proposition}\label{isogen} Let $\gamma: \widetilde G\to
G$ be an isogeny of nontrivial connected reductive
groups.~Then
\begin{equation*}{\rm gtd}(\widetilde
G)= {\rm gtd}(G).\end{equation*}
\end{proposition}
\begin{proof} If $G$ is a torus, then
$\widetilde G$ is a torus as well and the claim
follows from Proposition~\ref{1-2}(i). Assume now that
$G$ is nonabelian. By Proposition~\ref{1-2}(ii), there
is a parabolic subgroup $P$ of $\,\widetilde G$ such
that
\begin{equation}\label{=}
 {\rm gtd}(\widetilde G)={\rm
gtd}(\widetilde\alpha),\end{equation} where
 $\widetilde\alpha$ is the
 action of
$\,\widetilde G$ on ${\widetilde G}/P$. As $\gamma$ is
an isogeny, $\ker\gamma$ lies in the center of
$\widetilde G$ that in turn lies in $P$,
\cite[7.6.4]{springer}. So $\ker\gamma$ acts trivially
on $\widetilde G/P$ and hence $\widetilde\alpha$
descends to the action $\alpha$ of $\,G$ on
$\widetilde G/P$ such that ${\rm gtd}(\alpha)={\rm
gtd}(\widetilde\alpha)$. The claim now follows from
\eqref{=}, Definition~\ref{defgtdeg}, and
Lemma~\ref{>}. \quad $\square$
\renewcommand{\qed}{}
\end{proof}

Since every nonabelian connected reductive group $G$
admits an isogeny
$$Z\times S_1\times\ldots\times S_d\to
G,$$ where $Z$ is a torus and each $S_i$ is a simply
connected simple algebraic group, \cite[22.9,
22.10]{borel}, \cite[8.1.11, 10.1.1]{springer},
Propositions~\ref{1-2} and~\ref{isogen} reduce
calculating generic transitivity degrees to connecting
reductive groups to calculating the numbers ${\rm
gtd}(Z\times S_1\times\ldots\times S_d)$.

\begin{proposition}\label{product}
In the previous notation, there is an index $j$ and a
proper maximal parabolic subgroup $P_j$ of $S_j$ such
that
$${\rm gtd}(Z\times
S_1\times\ldots\times S_d)={\rm gtd}(\beta_j),$$ where
$\beta_j$ is the action of $Z\times
S_1\times\ldots\times S_d$ on $S_j/P_j$ given by $(z,
s_1,\ldots,s_d)\cdot x:= s_j\cdot x$.
\end{proposition}
\begin{proof}
As all proper  maximal parabolic subgroups of $Z\times
S_1\times\ldots\times S_d$ are exactly the subgroups
obtained from $Z\times S_1\times\ldots\times S_d$ by
replacing some $S_i$ with a proper maximal parabolic
subgroup  of $S_i$ (see Section~\ref{stpar}), the
claim follows from Proposition~\ref{1-2}(ii).
 \quad $\square$
\renewcommand{\qed}{}
\end{proof}

\noindent{\bf Corollary.} {\it In the previous
notation},
\begin{equation*}
{\rm gtd}(Z\times S_1\times\ldots\times
S_d)=\underset{i}{\rm max}\;{\rm gtd}(S_i).\quad
\square
\end{equation*}

\section{\bf Standard parabolic
subgroups}\label{stpar}

In this section we collect some necessary known facts
about parabolic subgroups.

 Let $G$ be a connected
reductive group. Fix a maximal torus $T$ of $\,G$. Let
$\Phi$, $\Phi^+$ and
$\Delta=\{\alpha_1,\ldots,\alpha_r\}$ be respectively
the root system of $\,G$ with respect to $T$, the
system of positive roots and the system of simple
roots of $\Phi$ determined by a fixed Borel subgroup
containing $T$. For a root $\alpha\in \Phi$, let
$u_\alpha$ be the one-dimensional unipotent root
subgroup of $\,G$ corresponding to $\alpha$.

If $I$ is a subset of $\Delta$, denote by $\Phi_I$ the
set of roots that are linear combinations of the roots
in $I$. Let $L_I$ be the subgroup of $\,G$ generated
by $T$ and all the $u_\alpha$'s with $\alpha\in
\Phi_I$. Let $U_I$ (respectively, $U_I^{-}$) be the
subgroup of $\,G$ generated by all $u_{\alpha}$ with
$\alpha\in \Phi^+\setminus \Phi_I$ (respectively,
$-\alpha\in \Phi^+\setminus \Phi_I$).
 Then $P_I:=L_IU_I$ and
 $P_I^-:=L_IU_I^-$ are
 parabolic subgroups of $\,G$ opposite
 to one another, $U_I$ and $U_I^-$
are the unipotent radicals of $P_I$ and $P_I^-$
respectively, $L_I$ is a Levi subgroup of $P_I$ and
$P_I^-$. In particular,
\begin{equation*}\label{udim}
\dim G= \dim L^{}_I+2\dim U_I^-.
\end{equation*}

Every parabolic subgroup of $\,G$ is conjugate to a
unique $P_I$, called {\it standard},
\cite[8.4.3]{springer}. We denote by $\mathfrak u_I^-$
the Lie algebra of $U_I^-$. For $I=\Delta\setminus
\{\alpha_i\}$, we denote $L_I$, $U_I^-$, $\mathfrak
u_I^-$, and $P_I$ respectively by $L_i$, $U_i^-$,
$\mathfrak u_i^-$, and $P_i$. Up to conjugacy,
$P_1,\ldots, P_r$ are all nonconjugate proper  maximal
parabolic subgroups of $\,G$.

Let $w_0$ be the element of the Weyl group of $G$ with
maximal length (regarding $\Delta$), and let
$\overset{.}w_0$ be a representative of $w_0$ in the
normalizer of $T$. Then there is an automorphism
$\varepsilon$ of $\Phi$ such that
$\varepsilon(\Delta)=\Delta$ and $\overset{.}w_0
P_I^-\overset{.}w_0^{-1} =P_{\varepsilon(I)}$. If $G$
is simple, then $\varepsilon$ is given by Table 6,
cf.\;\cite{bourbaki}:

\hskip 1mm

\begin{center}
\begin{tabular}{c|c}
\text{\rm type of} $G$ & $\varepsilon$
\\[2pt]
\hline \hline
&\\[-4pt]
${\sf A}_l$& $\begin{matrix}
 \varepsilon(\alpha_i)=\alpha_{l+1-i}
 \end{matrix}$
 \mbox{\ {\fontsize{9pt}{5mm} \selectfont\rm for all}
$i$}
\\[6pt]
 \hline
&\\[-6pt]
 ${\sf B}_l,$ ${\sf C}_l$, ${\sf E}_7$,
 ${\sf E}_8$, ${\sf F}_4$, ${\sf G}_2$
 & ${\rm id}$
\\[8pt]
\hline
&\\[-8pt]
${\sf D}_l$ & $\begin{matrix}
 {\rm id}
 \mbox{\
{\fontsize{9pt}{5mm} \selectfont\rm if}
$l$
{\fontsize{9pt}{5mm} \selectfont\rm is even,}}\\[1pt]
\begin{matrix}
\varepsilon(\alpha_{l-1})=\alpha_l,\
\varepsilon(\alpha_{l})=\alpha_{l-1},\\
\varepsilon(\alpha_i)=\alpha_i \mbox{
{\fontsize{9pt}{5mm} \selectfont\rm for} $i\neq l-1,
l$}
\end{matrix}\Bigr\}
\mbox{\ {\fontsize{9pt}{5mm} \selectfont\rm if} $l$
{\fontsize{9pt}{5mm} \selectfont\rm is odd}}
\end{matrix}$
\\[23pt]
\hline
&\\[-8pt]
${\sf E}_6$& $\begin{matrix}
\varepsilon(\alpha_1)=\alpha_6,\
\varepsilon(\alpha_2)=\alpha_2,\
\varepsilon(\alpha_3)=\alpha_5,\\
\varepsilon(\alpha_4)=\alpha_4,\
\varepsilon(\alpha_5)=\alpha_3,\
\varepsilon(\alpha_6)=\alpha_1
\end{matrix}$
\\
\end{tabular}
\nopagebreak \vskip 6mm
 \centerline{\rm Table $6$}
\end{center}

\vskip 2mm

So $P_I$ is conjugate to $P_I^-$ if and only if
$\varepsilon(I)=I$. In particular, $P_i$ is conjugate
to $P_i^-$ if and only if
$\varepsilon(\alpha_i)=\alpha_i$.

By the argument from the proof of
Proposition~\ref{isogen}, replacing $G$ by an
isogenous group does not change ${\rm gtd}(G:G/P_I)$.
On the other hand, if $G$ is a product of simply
connected simple algebraic groups and a torus, then
$\varepsilon $ is induced by an automorphism of $G$
stabilizing $T$, cf.\;\cite[32.1]{humphreys}. Hence,
by Lemma~\ref{>}(i), for every $G$ and $I$, we have
\begin{equation}\label{pp}
{\rm gtd}(G:P_I)={\rm gtd}(G:P_I^-).
\end{equation}

The center of reductive group $L_I$ is
$(r-|I|)$-dimensional. The root system of $L_I$ with
respect to
 $T$ is $\Phi_I$. The
subgroup of $L_I$ generated by $T$ and all the
$u_{\alpha}$\hskip -.5mm's with $\alpha\in
\Phi_I^+:=\Phi^+\cap \Phi_I$ is a Borel subgroup $B_I$
of $L_I$.
 The
systems of positive roots and simple roots of $\Phi_I$
determined by $B_I$ are respectively $\Phi_I^+$ and
$I$. Thus the Dynkin diagram of $(L_I, L_I)$ is
obtained from that of $(G, G)$ by removing the nodes
corresponding to the elements of $\Delta\setminus I$
togeher with the adjacent edges. In the sequel,
``roots'', ``positive roots'' and ``simple roots'' of
$L_I$ mean the elements of $\Phi_I$, $\Phi_I^+$ and
$I$ respectively, and highest weights of irreducible
$L_I$-modules are taken with respect to $T$ and~$B_I$.

The $L$-module structure of $\mathfrak u_I^-$
determined by the adjoint action is described as
follows, \cite{abs} (see also \cite{roh}). Every such
module is a direct sum of pairwise nonisomorphic
irreducible submo\-du\-les. Every such submodule is
completely (up to isomorphism) determined by its
highest weight. To describe these highest weights, it
is convenient to use the following terminology and
notation. Given a root $\beta\in \Phi^+\setminus
\Phi_I^+$, write it as a linear combination of simple
roots,
$$\textstyle \beta=\sum_{\alpha\in
I}a_{\alpha}\alpha+ \sum_{\alpha\in \Delta\setminus
I}b_{\alpha}\alpha.$$ Call $\sum_{\alpha\in
I}a_{\alpha} + \sum_{\alpha\in \Delta\setminus
I}b_{\alpha}$ the {\it height}, $\sum_{\alpha\in
\Delta\setminus I}b_{\alpha}\alpha$ the {\it shape},
and $\sum_{\alpha\in \Delta\setminus I}b_{\alpha}$ the
{\it level} of~$\beta$. Among all roots $\beta$ having
the same shape there is a unique one, $\beta_0$, whose
height is minimal.  Then $-\beta_0$ is the highest
weight of one of the irreducible submodules of the
$L_I$-module $\mathfrak u_I^-$; the shape and level of
$\beta_0$ are called the shape and level of this
submodule. The action of the center of $L_I$ on this
submodule is nontrivial. The highest weight of every
irreducible submodule $M$ of the $L_I$-module
$\mathfrak u_I^-$  is obtained in this way. The sum of
all the $M$'s of level $i$ is isomorphic to the
$L_I$-module $(\mathfrak u_I^-)^{(i)}/( \mathfrak
u_I^-)^{(i+1)}$, where $(\mathfrak u_I^-)^{(i)}$ is
the $i$-th term of the lower central series of
$\mathfrak u_I^-$. By \cite{richardson}, there are
only finitely many $L_I$-orbits in each $(\mathfrak
u_I^-)^{(i)}/( \mathfrak u_I^-)^{(i+1)}$.

According to this description, for $L_i$, only the
shapes of the form $b\alpha_i$ may occur, so
$\alpha_i$ is the unique root of level $1$ and
$\mathfrak u_I^-/[ \mathfrak u_I^-, \mathfrak u_I^-]$
is an irreducible $L_i$-module with the highest
weight $-\alpha_i$. In particular, if $\mathfrak
u_I^-$ is abelian, $\mathfrak u_I^-$ is an irreducible
$L_i$-module with the highest weight~$-\alpha_i$.

\smallskip

In Sections \ref{Al}--\ref{G2} we shall find ${\rm
gtd}(L_i:\mathfrak u_i^-)$ for every connected simple
algebraic group $G$ and every $i=1,\ldots, r$.

\section
{\bf \boldmath ${\rm gtd}(L_i:\mathfrak  u_i^-)$ for
\boldmath$G$ of type \boldmath${\sf A}_l$}\label{Al}

\begin{proposition}\label{al}
Let $\,G$ be a connected simple algebraic group of
type ${\sf A}_l$. Then
\begin{equation*}
{\rm gtd}(L_i:\mathfrak u_i^-)=
\begin{cases}1&
\text{if \ $2i=l+1$},\\
s_{li}&\text{if \ $2i\neq l+1$}
\end{cases}
\end{equation*}
where $s_{li}$ is defined by formulas \eqref{sli},
\eqref{Sli}.
\end{proposition}
\begin{proof}
We may (and shall) assume that
\begin{gather}\label{slp}
G={\bf SL}_{l+1}, \hskip 3mm P_i=\bigl\{
\begin{bmatrix}
A&B\\0&C
\end{bmatrix}
\in {\bf SL}_{l+1} \mid A\in {\bf GL}_i\big\},\\
L_i=\bigl\{
\begin{bmatrix}
A&0\\0&C
\end{bmatrix}
\in {\bf SL}_{l+1} \mid A\in {\bf GL}_i\bigr\}, \hskip
3mm \mathfrak u_i^-={\rm Mat}_{(l+1-i)\times i}.\notag
\end{gather}
The action of $L_i$ on $\mathfrak u_i^-$ is given by
\begin{equation*}
\begin{bmatrix}
A&0\\0&C\end{bmatrix} \cdot X:=
CXA^{-1}.\end{equation*} This yields
\begin{equation}\label{dimensions}
\dim L_i=2i^2+l^2-2li+2l-2i, \qquad \dim \mathfrak
u_i^-=il-i^2+i
\end{equation}
and shows that the action of $L_i$ on $(\mathfrak
u_i^-)^{\oplus a}$ (the direct sum
 of $a$ copies of
$\mathfrak u_i^-$) is equivalent in the sense of
\cite[Definition 4, p.\,36]{sk} to the action of ${\bf
GL}_i\times{\bf GL}_{l+1-i}$ on the space of $(i,
l+1-i)$-dimensional representations of the quiver

\vskip -1mm

\begin{equation*}
1\hskip .5mm\begin{matrix} \xymatrix@=16mm@M=
-.4mm{\circ \ar@/^.8pc/[r]_{\hskip .5mm \colon}
\ar@{->}[r] \ar@/_.8pc/[r]^{\hskip .5mm \colon} &\circ
}
\end{matrix}\hskip .5mm 2
\hskip 1mm
\end{equation*}

\

\noindent with $a$ arrows. From \cite[Theorem~4]{kac1}
we obtain that for $a\geqslant 2$ an open orbit in
this space exists if and only if (a) ${\mathcal
S}_{li}\neq \varnothing$ and (b) $a\in {\mathcal
S}_{li}$. Since condition (a) is equivalent to the
inequality $2i\neq l+1$, we now deduce from
\eqref{sli} and Definition~\ref{defgtdeg}
 that ${\rm gtd}(L_i:\mathfrak
u_i^-)=s_{li}$ if $2i\neq l+1$.  If $2i=l+1$, then
\eqref{dimensions} yields $\dim L_i<2\dim \mathfrak
u_i^-$, and hence in this case there is no open
$L_i$-orbit in $(\mathfrak u_i^-)^{\oplus a}$ for
$a\geqslant 2$. On the other hand, if $a=1$, then, by
Remark~\ref{abel}, such an orbit exists. \quad
$\square$
\renewcommand{\qed}{}\end{proof}

\section
{\bf \boldmath ${\rm gtd}(L_i:\mathfrak  u_i^-)$ for
\boldmath$G$ of type \boldmath${\sf B}_l$, $l\geqslant
3$ }\label{gbl}

\begin{proposition}\label{bl} Let $G$ be
a connected simple algebraic group of type ${\sf
B}_l$, $l\geqslant 3$. Then
\begin{equation*}
{\rm gtd}(L_i:\mathfrak u_i^-)=\begin{cases}0& \mbox
{if}\hskip 2mm i\neq 1,l,\\ 1&\mbox{if}\hskip 2mm i=1,
l.\end{cases}
\end{equation*}
\end{proposition}
\begin{proof}
{\it Step} 1.\;Let $i=1$. By Remark \ref{abel} and
Corollary~\ref{abelianrad} of Proposition~\ref{p-},
the action of $L_1$ on $\mathfrak  u_1^-$ is locally
transitive. The type of $(L_1, L_1)$ is ${\sf
B}_{l-1}$. The action of $L_1$ on $\mathfrak u_1^-$ is
irreducible with the highest weight $-\alpha_1=(-2,
1,0,\ldots,0)$. Hence $\mathfrak u_1^-$, considered as
(irreducible) $(L_1, L_1)$-module,
 has the highest
weight with the Dynkin diagram
\begin{equation}\label{diag-bl1}
\xymatrix@=4.5mm@M= -.4mm{\circ\ar@{-}[r]^{\hskip -
8.1mm 1} &{\ \ldots\ } \ar@{-}[r]&\circ\ar@{=>}[r]&
\circ}\quad .
\end{equation}

From \eqref{diag-bl1} we deduce that, for a general
point $z\in \mathfrak u_1^-$, the group $(L_1,
L_1)_z^0$ is locally isomorphic to ${\bf SO}_{2l-1}$,
and the codimension of $(L_1, L_1)\cdot z$ in
$\mathfrak u_1^-$ is equal to $1$. As the
(one-dimensional) center of $L_1$ acts on $\mathfrak
u_1^-$ nontrivially and the action of $L_1$ on
$\mathfrak u_1^-$ is locally transitive, this yields
$(L_1)_z^0=(L_1, L_1)_z^0$. Hence the action of
$(L_1)_z^0$ on $\mathfrak u_1^-$ is not locally
transitive. Lemma~\ref{reduction} now yields ${\rm
gtd}(L_1:\mathfrak u_1^-)=1$.

\smallskip

{\it Step} 2.\;Let $i=l$. The type of $(L_l, L_l)$ is
${\sf A}_{l-1}$. By inspection  of $\Phi^+$ in
\cite[Planche II]{bourbaki} we obtain  that,  for the
action of $L_l$ on $\mathfrak u_l^-$, there are
exactly two shapes $\alpha_l$ and $2\alpha_l$, and
they determine the highest weights $-\alpha_l=
(0,\ldots, 0,1,-2)$ and $-\alpha_{l-1}-2\alpha_l=
(0,\ldots,0, 1, 0, -2)$. Hence $\mathfrak u_1^-$ is
the direct sum of two irreducible $L_l$-modules
$\mathfrak u_{l1}^-$ and $\mathfrak u_{l2}^-$ that,
considered as (irreducible) $(L_l, L_l)$-modules, have
respectively the highest weights with the Dynkin
diagrams
\begin{gather}\label{diag-bl12}
\begin{gathered}
\xymatrix@=4.5mm@M= -.4mm{\circ\ar@{-}[r] &{\ \ldots\
} \ar@{-}[r]^{\hskip 8.1mm 1}&\circ} \hskip 6mm
\text{and} \hskip 6mm \xymatrix@=4.5mm@M=
-.4mm{\circ\ar@{-}[r] &{\ \ldots\ } \ar@{-}[r]^{\hskip
8.1mm 1}&\circ\ar@{-}[r]&\circ }\quad .
\end{gathered}
\end{gather}

\vskip 2mm

From \eqref{diag-bl12} and \cite[Section 3, A(8)(iii)
and B(4)(iii)]{kimu1} we conclude that the action of
$L_l$ on $\mathfrak u_l^-$ is locally transitive. On
the other hand, if the action of $L_l$ on $\mathfrak
u_l^- \oplus \mathfrak u_l^-$ would be locally
transitive, then all the more the action of ${\bf
GL}_1^4\times (L_l, L_l)$ on $\mathfrak u_{l1}^-\oplus
\mathfrak u_{l1}^-\oplus \mathfrak u_{l2}^-\oplus
\mathfrak u_{l2}^-$, where ${\bf GL}_1^4$ acts on the
direct summands by independent scalar multiplications,
would be locally transitive. Using the terminology and
notation of \cite{kimu1}, this would mean that the
pair $({\bf GL}_1^4\times {\bf SL}_l, \Lambda_1\oplus
\Lambda_1\oplus\Lambda_2\oplus\Lambda_2)$ is
prehomogeneous. However the classification obtained in
\cite[Section~3]{kimu1} shows that it is not so. Thus
${\rm gtd}(L_l:\mathfrak u_l^-)=1$.

\smallskip

{\it Step} 3. The type of $(L_i, L_i)$ for $1<i<l$ is
${\sf A}_{i-1}+{\sf B}_{l-i}$, where ${\sf B}_1:={\sf
A}_1$. By inspection of $\Phi^+$ in \cite[Planche
II]{bourbaki}  we obtain that  for the action of $L_i$
on $\mathfrak u_i^-$ there are exactly two shapes
$\alpha_i$ and $2\alpha_i$, and they determine
respectively the highest weights
\begin{equation*}-\alpha_i=\begin{cases}
(0,\ldots,0, 1, -\underset{ \hat i}{2}, 1,0,\ldots,
0)& \text{
if $i\neq l-1$,}\\
(0,\ldots,0,1,-2,2)&\text{ if $i=l-1$,}
\end{cases}\quad \text{and}
\end{equation*}
\begin{equation*}
-\alpha_{i-1}-2\alpha_{i}-\ldots-2\alpha_l
=\begin{cases} (0,\ldots,0,1,0, -\underset{ \hat
i}{1},0,\ldots,0)&\text{if $i\neq 2$},\\
 (0,-1,0,\ldots,0)&\text{if $i=2$}.
 \end{cases}
 \end{equation*}
 Hence
$\mathfrak u_i^-$ is the direct sum  of two
irreducible $L_i$-modules $\mathfrak u_{i1}^-$ and
$\mathfrak u_{i2}^-$ that, considered as (irreducible)
$(L_i, L_i)$-modules, have respectively the highest
weights with the Dynkin diagrams
\begin{gather}\label{cli>1}
\begin{gathered}
\xymatrix@=4.5mm@M= -.4mm{\circ\ar@{-}[r]&{\ \ldots\
}\ar@{-}[r]^{\hskip 8mm 1}&\circ&\circ
\ar@{-}[r]^{\hskip -8.5mm 1}& {\ \ldots \
}\ar@{-}[r]&\circ\ar@{=>}[r]&\circ}\hskip 5mm
\text{and} \hskip 5mm \xymatrix@=4.5mm@M=
-.4mm{\circ\ar@{-}[r]&{\ \ldots\ }\ar@{-}[r]^{\hskip
8mm 1} &\circ\ar@{-}[r]&\circ&\circ \ar@{-}[r]& {\
\ldots \ }\ar@{-}[r]&\circ\ar@{=>}[r]&\circ}\hskip 4mm
\text{if
$i\neq 2, l-1$,}\\[-2pt]
\xymatrix@=4.5mm@M= -.4mm{\circ\ar@{-}[r]&{\ \ldots\
}\ar@{-}[r]^{\hskip 8mm 1}&\circ&\circ
\ar@{}[r]^{\hskip -4.5mm 2}& }\hskip 1mm \text{and}
\hskip 4.2mm \xymatrix@=4.5mm@M=
-.4mm{\circ\ar@{-}[r]&{\ \ldots\ }\ar@{-}[r]^{\hskip
8mm 1} &\circ\ar@{-}[r]&\circ&\circ \ar@{}[r]& }\hskip
1mm \text{if $i=l-1$, $l>3$,}\\[-2pt]
\xymatrix@=4.5mm@M= -.4mm{\circ\ar@{}[r]^{\hskip
-5.8mm 1} &\circ\ar@{-}[r]^{\hskip -8.5mm 1}& {\
\ldots \ }\ar@{-}[r]&\circ\ar@{=>}[r]&\circ}\hskip 5mm
\text{and} \hskip 5mm \xymatrix@=4.5mm@M= -.4mm{\circ
&\circ\ar@{-}[r]& {\ \ldots \
}\ar@{-}[r]&\circ\ar@{=>}[r]&\circ} \hskip 4mm
\text{if\quad $i=2$, $l> 3$,}
\\[-2pt]
\xymatrix@=4.5mm@M= -.4mm{\circ\ar@{}[r]^{\hskip
-5.7mm 1} &\circ\ar@{}[]^{\hskip -.1mm 2}&} \hskip 4mm
\text{and} \hskip 4mm \xymatrix@=4.5mm@M=
-.4mm{\circ&\circ} \hskip 4mm \text{if $i=2$, $l=3$.}
\end{gathered}
\end{gather}

\smallskip

{\it Step} 4.\;Let $i=2$. From \eqref{cli>1} we deduce
that $\dim\mathfrak u_{22}^-=1$. Since the action of
$L_2$ on $\mathfrak u_{22}^-$ is locally transitive
and the center of $L_2$ is one-dimensional,
$(L_2)_z^0=(L_2, L_2)$ for a nonzero $z\in \mathfrak
u_{22}^-$. This and \eqref{cli>1} now yield that, in
the notation and terminology of \cite{sk},
\cite{kkiy}, \cite{kkti}, the action of $(L_2)_z^0$ on
$\mathfrak u_{21}^-$ is equivalent to that determined
by the pair $({\bf SL}_2\times {\bf SO}_{2l-3},
\Lambda_1\otimes\Lambda_1)$. By \cite[\S\,6]{sk}, this
pair is not prehomogeneous. So the action of
$(L_2)_z^0$ on $\mathfrak u_{21}^-$ is not locally
transitive. Hence, by Lemma~\ref{reduction}, the
action of $L_2$ on $\mathfrak u_{2}^-$ is not locally
transitive as well.

\smallskip

{\it Step} 5.\;Let $i=l-1$ and $i\neq 2$ (hence
$l>3$). From \eqref{cli>1} we deduce that the action
of $L_{l-1}$ on $\mathfrak u_{l-1}^-$ is equivalent to
that determined by the pair $({\bf GL}_{l-1}\times
{\bf SO}_3, \Lambda_1\otimes\Lambda_1+\Lambda_2\otimes
1)$. This shows that the action of $L_{l-1}$ on
$\mathfrak u_{l-1, 2}^-$ is equivalent to the natural
action of ${\bf GL}_{l-1}$ on the space of
skew-symmetric bilinear forms over $k$ in $l-1$
variables.

Assume that $l$ is odd. The last remark then yields
that, for a general point $z\in \mathfrak u_{l-1,
2}^-$, the group $(L_{l-1})_z^0$ is locally isomorphic
to ${\bf Sp}_{l-1}\times {\bf SO}_3$ and its action on
$\mathfrak u_{l-1, 1}^-$ is equivalent to that
determined by the pair $({\bf Sp}_{l-1}\times {\bf
SO}_3, \Lambda_1\otimes\Lambda_1)$. By
\cite[\S\,6]{sk}, this pair is not
pre\-ho\-mo\-ge\-neous. So the action of
$(L_{l-1})_z^0$ on $\mathfrak u_{l-1,1}^-$ is not
locally transitive. Hence, by Lemma~\ref{reduction},
if $l$ is odd, then the action of $L_{l-1}$ on
$\mathfrak u_{l-1}^-$ is not locally transitive.

Assume now that $l$ is even. Then, according to
\cite[Proposition\,1.34\,(2)]{kkti}, the
pre\-ho\-mo\-ge\-neity of $({\bf GL}_{l-1}\times {\bf
SO}_3, \Lambda_1\otimes\Lambda_1+\Lambda_2\otimes 1)$
is equivalent to that of $({\bf GL}_{2}\times {\bf
SO}_3,\linebreak \Lambda_1\otimes\Lambda_1+{\rm
det}\otimes 1)$. In turn, since for $({\bf
GL}_{2}\times {\bf SO}_3, {\rm det}\otimes 1)$ the
stabilizer of a general point is ${\bf SL}_{2}\times
{\bf SO}_3$, Lemma~\ref{reduction} shows that the
prehomoge\-neity of $({\bf GL}_{2}\times {\bf SO}_3,
\Lambda_1\otimes\Lambda_1+{\rm det}\otimes 1)$ is
equivalent to that of $({\bf SL}_{2}\times {\bf SO}_3,
\Lambda_1\otimes\Lambda_1)$. But according to
\cite[\S\,6]{sk},  the last pair is not
prehomogeneous.

Summing up, we obtain that for every $l$ the action of
$L_{l-1}$ on $\mathfrak u_{l-1}^-$ is not locally
transitive.

\smallskip

{\it Step} 6.\;Let $2<i<l-1$. From \eqref{cli>1} we
deduce that the action of $L_i$ on $\mathfrak u_i^-$
is equivalent to that determined by the pair $( {\bf
GL}_i\times {\bf SO}_{2(l-i)+1},
\Lambda_1\otimes\Lambda_1+\Lambda_2\otimes 1)$. If the
action of $L_i$ on $\mathfrak u_i^-$ would be locally
transitive, then all the more the action determined by
the pair $({\bf GL}_1^2\times {\bf SL}_i\times {\bf
SO}_{2(l-i)+1},
\Lambda_1\otimes\Lambda_1+\Lambda_2\otimes 1)$, where
${\bf GL}_1^2$ acts on the summands by independent
scalar multiplications, would be locally transitive.
In turn, this would mean that if $2(l-i)+1>i$, then
the last pair is $2$-simple prehomogeneous of type I
in the sense of \cite{kkiy}. However the
classification of such pairs obtained in \cite[Section
3]{kkiy} shows that is is not so. Hence, for
$2(l-i)+1>i$, the action of $L_i$ on $\mathfrak u_i^-$
is not locally transitive.

Assume now that $2(l-i)+1\leqslant i$. Then, by the
same argument, if the action of $L_i$ on $\mathfrak
u_i^-$ would be locally transitive, then the pair
$({\bf GL}_1^2\times {\bf SL}_i\times {\bf
SO}_{2(l-i)+1},
\Lambda_1\otimes\Lambda_1+\Lambda_2\otimes 1)$ would
be $2$-simple prehomogeneous
 of type II in the sense of
\cite{kkiy}, \cite{kkti}. However the classification
of such pairs obtained in \cite[Section 5]{kkti} shows
that it is not so. Hence if $2(l-i)+1\leqslant i$,
then the action of $L_i$ on $\mathfrak u_i^-$ is not
locally transitive.

Summing up, we obtain that if $2<i<l-2$, then the
action of $L_i$ on $\mathfrak u_i^-$ is not locally
transitive. \quad $\square$
\renewcommand{\qed}{}\end{proof}

\begin{remark} In \cite{kime1},
\cite{kime2}, for every $i$, it is obtained a
classification of all connected reductive subgroups of
$G={\bf SO}_n(\mathbb C)$ that act locally
transitively on $G/P_i$. \quad $\square$ \end{remark}


\section{\bf \bf
\boldmath ${\rm gtd}(L_i:\mathfrak  u_i^-)$
for \boldmath$G$ of type \boldmath${\sf C}_l$,
$l\geqslant 2$}\label{cl>1}

\begin{proposition}\label{cl} Let $G$ be
a connected simple algebraic group of type ${\sf
C}_l$, $l\geqslant 2$. Then
\begin{equation*}
{\rm gtd}(L_i:\mathfrak u_i^-)=\begin{cases}0& \mbox
{if}\hskip 2mm i\neq 1,l,\\ 1&\mbox{if}\hskip 2mm i=1,
l.\end{cases}
\end{equation*}
\end{proposition}
\begin{proof}
{\it Step} 1.\;Let $i=1$. The type of $(L_1, L_1)$ is
${\sf C}_{l-1}$, where ${\sf C}_1:={\sf A}_1$. By
inspection of $\Phi^+$ in \cite[Planche III]{bourbaki}
we obtain that,  for the action of $L_1$ on $\mathfrak
u_1^-$, there are exactly two shapes $\alpha_1$ and
$2\alpha_1$, and they determine the highest weights
$-\alpha_1= (-2, 1,0,\ldots,0)$ and
$-2\alpha_1-\ldots-2\alpha_{l-1}-\alpha_l=
(-2,0,\ldots,0)$. Hence $\mathfrak u_1^-$ is the
direct sum of two irreducible $L_1$-modules $\mathfrak
u_{11}^-$ and $\mathfrak u_{12}^-$, where $\mathfrak
u_{11}^-$, considered as (irreducible) $(L_1,
L_1)$-module,  has the highest weight with the Dynkin
diagram
\begin{gather}\label{diag-cl1}
\begin{gathered}
\xymatrix@=4.5mm@M= -.4mm{\circ\ar@{-}[r]^{\hskip -
8.1mm 1} &{\ \ldots\ } \ar@{-}[r]&\circ\ar@{<=}[r]&
\circ} \quad
\text{if $l>2$,}\\
\overset{1}{\circ}\quad \text{if $l=2$,}
\end{gathered}
\end{gather}
and $\mathfrak u_{12}^-$ is a trivial one-dimensional
module.

As the action of $L_1$ on $\mathfrak u_{12}^-$ is
locally transitive and $(L_1, L_1)$ has no nontrivial
characters, $(L_1)_z^0=(L_1, L_1)$ for a nonzero point
$z\in \mathfrak u_{12}^-$. It follows from
\eqref{diag-cl1} that the action of $(L_1, L_1)$ on
$\mathfrak u_{11}^-$ is equivalent, in the sense of
\cite[Definition 4, p.\,36]{sk}, to the natural action
of ${\bf Sp}_{2l-2}$ on $k^{2l-2}$. As the latter is
locally transitive by Witt's theorem,
Lemma~\ref{reduction} yields that the action of $L_1$
on $\mathfrak u_1^-$ is locally transitive. Since
${\bf Sp}_{2l-2}$ fixes  a nondegenerate
skew-symmetric form on $k^{2l-2}$, the natural action
of ${\bf Sp}_{2l-2}$ on $k^{2l-2}\oplus k^{2l-2}$ is
not locally transitive. Applying the same argument as
above we then conclude that ${\rm gtd}(L_1:\mathfrak
u_1^-)=1$.

\smallskip

{\it Step} 2.\;Let $i=l$. By Remark \ref{abel} and
Corollary~\ref{abelianrad} of Proposition~\ref{p-},
the action of $L_l$ on $\mathfrak  u_l^-$ is locally
transitive.  The type of $(L_l, L_l)$ is ${\sf
A}_{l-1}$, so $\dim L_l=l^2$, $\dim \mathfrak
u_l^-=l(l+1)/2$. As $\dim L_l<2\dim \mathfrak u_l^-$,
 we have ${\rm gtd}(L_l:\mathfrak u_l^-)=1$.

\smallskip

{\it Step} 3.\;Let $l\geqslant 3$ and $1<i<l$. The
type of $(L_i, L_i)$ is ${\sf A}_{i-1}+{\sf C}_{l-i}$.
By inspection of $\Phi^+$ in \cite[Planche
III]{bourbaki} we obtain that,  for the action of
$L_i$ on $\mathfrak u_i^-$, there are exactly two
shapes $\alpha_i$ and $2\alpha_i$, and they determine
respectively the highest weights $-\alpha_i=
(0,\ldots,0, 1, -\underset{ \hat i}{2}, 1,0,\ldots, 0)
$ and $-2\alpha_i-\ldots-2\alpha_{l-1}-\alpha_l=
(0,\ldots,0,2,-\underset{ \hat i}{2},0\ldots,0) $.
Hence $\mathfrak u_i^-$ is the direct sum  of two
irreducible $L_i$-modules $\mathfrak u_{i1}^-$ and
$\mathfrak u_{i2}^-$ that, considered as (irreducible)
$(L_i, L_i)$-modules, have respectively the highest
weights with the Dynkin diagrams
\begin{gather}\label{cl>2}
\begin{gathered}
\xymatrix@=4.5mm@M= -.4mm{\circ\ar@{-}[r]&{\ \ldots\
}\ar@{-}[r]^{\hskip 8mm 1}&\circ&\circ
\ar@{-}[r]^{\hskip -8.5mm 1}& {\ \ldots \
}\ar@{-}[r]&\circ\ar@{<=}[r]&\circ}\hskip 6mm
\text{and} \hskip 7mm \xymatrix@=4.5mm@M=
-.4mm{\circ\ar@{-}[r]&{\ \ldots\ }\ar@{-}[r]^{\hskip
8mm 2}&\circ&\circ \ar@{-}[r]& {\ \ldots \
}\ar@{-}[r]&\circ\ar@{<=}
[r]&\circ}\hskip 5mm \text{if $i\neq l-1$,}\\
\xymatrix@=4.5mm@M= -.4mm{\circ\ar@{-}[r]&{\ \ldots\
}\ar@{-}[r]^{\hskip 8mm 1}&\circ&\circ
\ar@{}[r]^{\hskip -4.6mm 1}& }\hskip 3mm \text{and}
\hskip 6mm \xymatrix@=4.5mm@M=
-.4mm{\circ\ar@{-}[r]&{\ \ldots\ }\ar@{-}[r]^{\hskip
8mm 2}&\circ&\circ} \hskip 5mm\text{if $i=l-1$,}
\end{gathered}
\end{gather}

From \eqref{cl>2} we deduce that $\dim \mathfrak
u_{i2}^-=i(i+1)/2$, for a general point $z\in
\mathfrak u_{i2}^-$ the group
 $(L_i, L_i)_z^0$
  is locally isomorphic to ${\bf
SO}_{i}\times {\bf Sp}_{2l-2i}$, and the codimension
of orbit $(L_i, L_i)\cdot z$ in $\mathfrak u_{i2}^-$
is equal to $1$. As the center of $L_i$ is
one-dimensional  and acts on $\mathfrak u_{i2}^-$
nontrivially, and the action of $L_i$ on $\mathfrak
u_{i2}^-$ is locally transitive, this yields
$(L_i)_z^0=(L_i, L_i)_z^0$. So $(L_i)_z^0$ is locally
isomorphic to ${\bf SO}_{i}\times {\bf Sp}_{2l-2i}$.
From \eqref{cl>2} we now deduce that, in the
terminology and notation of \cite[Definition 4,
p.\,36]{sk}, the action of $(L_i)_z^0$ on $\mathfrak
u_{i1}^-$ is equivalent to that determined by the pair
$({\bf SO}_{i}\times {\bf Sp}_{2l-2i},
\Lambda_1\otimes\Lambda_1)$. By \cite[\S\,7]{sk}, this
pair is not prehomogeneous. Lemma~\ref{reduction} now
yields that the action of $L_i$ on $\mathfrak u_i^-$
is not locally transitive.
 \quad $\square$
\renewcommand{\qed}{}\end{proof}

\section{\bf \bf \boldmath ${\rm gtd}(L_i:\mathfrak  u_i^-)$
for \boldmath$G$ of type \boldmath${\sf D}_l$,
$l\geqslant 4$} \label{Dl}

\begin{proposition}\label{dl} Let $G$ be
a connected simple algebraic group of type ${\sf
D}_l$, $l\geqslant 4$. Then
\begin{equation*}
{\rm gtd}(L_i:\mathfrak u_i^-)=\begin{cases}0& \mbox
{if}\hskip 2mm i\neq 1,l-1, l,\\ 1&\mbox{if}\
i=1,\\
1& \mbox {if}\hskip 2mm \mbox{$l$ is even and}\hskip 2mm i= l-1, l,\\
2& \mbox {if}\hskip 2mm  \mbox{$l$ is odd and}\hskip
2mm i= l-1, l.\end{cases}
\end{equation*}
\end{proposition}

\begin{proof}
{\it Step} 1.\; Let $i=1$. By Remark \ref{abel} and
Corollary~\ref{abelianrad} of Proposition~\ref{p-},
the action of $L_1$ on $\mathfrak u_1^-$ is locally
transitive. The type of $(L_1, L_1)$ is ${\sf
D}_{l-1}$, where ${\sf D}_3:={\sf A}_3$. The action of
$L_1$ on $\mathfrak u_1^-$ is irreducible with the
highest weight $-\alpha_1=(-2, 1,0,\ldots,0)$. Hence
$\mathfrak u_1^-$, considered as (irreducible) $(L_1,
L_1)$-module,
 has the highest
weight with the Dynkin diagram
\begin{equation}\label{dl1}
\begin{matrix}
\xymatrix@=4.5mm@M= -.4mm@R=3mm{ &&&\circ
\\
\circ\ar@{-}[r]^{\hskip -8.22mm 1} & {\ \ldots \
}\ar@{-}[r] & \circ\ar@{-}[ur] \ar@{-}[dr]
\\
&&&\circ }
\end{matrix}
\quad .
\end{equation}

\vskip 2mm

Arguing like in Subsection~\ref{cl>1} for $l\geqslant
3$, we deduce from \eqref{dl1} that, for a general
point $z\in \mathfrak u_i^-$, the group $(L_1)_z^0$
lies in $(L_1, L_1)$. As the action of  $(L_1, L_1)$
on $\mathfrak u_1^-$ is not locally transitive (it
fixes a nondegenerate quadratic form),
Lemma~\ref{reduction} now yields  ${\rm
gtd}(L_1:\mathfrak u_1^-)=1$.

\smallskip

{\it Step} 2.\;Let $i=l-1, l$. Again by Remark
\ref{abel} and Corollary~\ref{abelianrad} of
Proposition~\ref{p-}, the action of $L_l$ on
$\mathfrak u_l^-$ is locally transitive. The type of
$(L_l, L_l)$ is ${\sf A}_{l-1}$. The action of $L_l$
on $\mathfrak u_l^-$ is irreducible with the highest
weight $-\alpha_l=(0,\ldots, 0,1,0,-2)$. Hence
$\mathfrak u_l^-$, considered as (irreducible) $(L_l,
L_l)$-module, has the highest weight with the Dynkin
diagram
\begin{gather}\label{diag-dll}
\begin{gathered}
\xymatrix@=4.5mm@M= -.4mm{\circ\ar@{-}[r] &{\ \ldots\
} \ar@{-}[r]^{\hskip 8.1mm 1}&\circ\ar@{-}[r]&\circ
}\quad .
\end{gathered}
\end{gather}

From \eqref{diag-dll} we deduce that, in the
terminology and notation of \cite[Definition 4,
p.\,36]{sk}, the action of $(L_l, L_l)$ on $\mathfrak
u_l^-\oplus \mathfrak u_l^-$ is equivalent to that
determined by the pair $({\bf SL}_l, \Lambda_2\oplus
\Lambda_2)$. Since the center of $L_l$ is
one-dimensional, it now follows from \cite[Proposition
2.2 and Section~3,~B~(3)(iii)]{kimu1} that the action
of $L_l$ on $\mathfrak u_l^-\oplus \mathfrak u_l^-$ is
not locally transitive for even $l$, and is locally
transitive for odd $l$. On the other hand, as $\dim
L_l=l^2<3l(l-1)/2=\dim (\mathfrak u_l^-\oplus
\mathfrak u_l^-\oplus \mathfrak u_l^-)$, the action of
$L_l$ on $\mathfrak u_l^-\oplus \mathfrak u_l^-\oplus
\mathfrak u_l^-$ is not locally transitive. So we see
that ${\rm gtd}(L_l:\mathfrak u_l^-)$ is equal to $1$
if $l$ is even, and to $2$ if $l$ is odd.

By Proposition~\ref{p-} and \eqref{pp}, we have ${\rm
gtd}(L_{l-1}:\mathfrak u_{l-1}^-)={\rm
gtd}(L_{l}:\mathfrak u_{l}^-)$.

\smallskip

{\it Step} 3.\;If $1<i\leqslant l-2$, then the type of
$(L_i, L_i)$ is ${\sf A}_{i-1}+{\sf D}_{l-i}$, where
${\sf D}_2:={\sf A}_1+{\sf A}_1$. By inspection of
\cite[Planche IV]{bourbaki}  we obtain  that  for the
action of $L_i$ on $\mathfrak u_i^-$ there are exactly
two shapes $\alpha_i$ and $2\alpha_i$, and they
determine respectively the highest weights
\begin{equation*}
-\alpha_i=\begin{cases} (0,\ldots, 1, -\underset{ \hat
i}{2}, 1,0,\ldots,0)&
\text{if $i\neq l-2$,}\\
(0,\ldots, 0,1,-2,1,1)&\text{if $i=l-2$,}
\end{cases}\quad \text{and}
\end{equation*}
\begin{equation*}
-\alpha_{i-1}-2\alpha_i-
\ldots-2\alpha_{l-2}-\alpha_{l-1}- \alpha_l=
\begin{cases}
(0,\ldots,0,1,0, -\underset{ \hat i}1,
0,\ldots,0)&\text{if $i\neq 2$},\\
(0,-1,0,\ldots,0)&\text{if $i=2$}.
\end{cases}
\end{equation*}
Hence $\mathfrak u_i^-$ is the direct sum of two
irreducible $L_i$-modules $\mathfrak u_{i1}^-$ and
$\mathfrak u_{i2}^-$ that, considered as (irreducible)
$(L_i, L_i)$-modules, have respectively the highest
weights with the Dynkin diagrams
\begin{gather}\label{diagrdi>1}
\begin{gathered}
\begin{matrix}\xymatrix@=4.5mm@M= -.4mm@R=3mm{
&&&&&&\circ\\
\circ\ar@{-}[r]&{\ \ldots \ } \ar@{-}[r]^{\hskip 8mm
1}&\circ &\circ\ar@{-}[r]^{\hskip -8.3mm 1}&{\ \ldots
\ }\ar@{-}[r]&\circ\ar@{-}[ur]
\ar@{-}[dr]&\\
&&&&&&\circ} \end{matrix}\hskip 6mm \text{and} \hskip
6mm
\begin{matrix}
\xymatrix@=4.5mm@M= -.4mm@R=3mm{
&&&&&&&\circ\\
\circ\ar@{-}[r]&{\ \ldots \ }
\ar@{-}[r]&\circ\ar@{-}[r]^{\hskip -5.8mm 1}&\circ
&\circ\ar@{-}[r]&{\ \ldots \
}\ar@{-}[r]&\circ\ar@{-}[ur]
\ar@{-}[dr]&\\
&&&&&&&\circ}\end{matrix}\quad\text{if
$i\neq 2, l-2$,}\\
\begin{matrix}
\xymatrix@=4.5mm@M= -.4mm@R=3mm{\circ\ar@{-}[r]&{\
\ldots \ } \ar@{-}[r]^{\hskip 8mm
1}&\circ&\circ\ar@{}^{\hskip -0.2mm
1}&\circ\ar@{}^{\hskip -.1mm 1}&}
\end{matrix}
\hskip 4mm \text{and} \hskip 7mm
\begin{matrix}
\xymatrix@=4.5mm@M= -.4mm@R=3mm{\circ\ar@{-}[r]&{\
\ldots \ } \ar@{-}[r]^{\hskip 8mm
1}&\circ\ar@{-}[r]&\circ &\circ&\circ }
\end{matrix}\quad\text{if $i=l-2$, $l>
4$,}\\
\begin{matrix}\xymatrix@=4.5mm@M= -.4mm@R=3mm{
&&&&\circ\\
\circ\ar@{}[r]^{\hskip -5.7mm 1}
&\circ\ar@{-}[r]^{\hskip -8.3mm 1}&{\ \ldots \
}\ar@{-}[r]&\circ\ar@{-}[ur]
\ar@{-}[dr]&\\
&&&&\circ} \end{matrix}\hskip 6mm \text{and} \hskip
6mm
\begin{matrix}
\xymatrix@=4.5mm@M= -.4mm@R=3mm{
&&&&\circ\\
\circ &\circ\ar@{-}[r]&{\ \ldots \
}\ar@{-}[r]&\circ\ar@{-}[ur]
\ar@{-}[dr]&\\
&&&&\circ}\end{matrix}\quad\text{if
$i=2$, $l>4$,}\\
\begin{matrix}
\xymatrix@=4.5mm@M= -.4mm@R=3mm{
\circ\ar@{}[r]^{\hskip -5.7mm
1}&\circ\ar@{}[r]^{\hskip -5.7mm 1}
&\circ\ar@{}[r]^{\hskip -4.7mm 1}&}
\end{matrix}
\hskip 4mm \text{and} \hskip 7mm
\begin{matrix}
\xymatrix@=4.5mm@M= -.4mm@R=3mm{ \circ&\circ&\circ}
\end{matrix}
\quad\text{if $i=2$, $l=4$}.
\end{gathered}
\end{gather}

{\it Step} 4.\;Let $i=2$. From \eqref{diagrdi>1} we
deduce that $\dim\mathfrak u_{22}^-=1$ and the action
of $(L_2, L_2)$ on $\mathfrak u_{21}^-$ is equivalent
to that determined by the pair $({\bf SL}_2\times {\bf
SO}_{2(l-1)}, \Lambda_1\otimes\Lambda_1)$ for $l>4$,
and to the pair $({\bf SL}_2\times{\bf SL}_2\times{\bf
SL}_2, \Lambda_1\otimes\Lambda_1\otimes\Lambda_1)$ for
$l=4$. By \cite[\S\,6]{sk}, these pairs are not
prehomogeneous. Using now the same argument as in the
case $i=2$ in Subsection~\ref{gbl} we obtain that the
action of $L_2$ on $\mathfrak u_2^-$ is not locally
transitive.

\smallskip

{\it Step} 5.\;Let $i=l-2$ and $i\neq 2$ (hence
$l>4$). From \eqref{diagrdi>1} we deduce that the
action of $L_{l-2}$ on $\mathfrak u_{l-2}^-$ is
equivalent to that determined by the pair $({\bf
GL}_{l-2}\times{\bf SL}_2\times{\bf SL}_2,
\Lambda_1\otimes\Lambda_1\otimes\Lambda_1+
\Lambda_2\otimes 1\otimes 1)$. Hence the action of
$L_{l-2}$ on $\mathfrak u_{l-2,2}^-$ is equivalent to
the natural action of ${\bf GL}_{l-2}$ on the space of
skew-symmetric bilinear forms over $k$ in $l-2$
variables.

Assume that $l$ is even. Then the last remark yields
that, for a general point $z\in \mathfrak
u_{l-2,2}^-$, the group $(L_{l-2})_z^0$ is locally
isomorphic to ${\bf Sp}_{l-2}\times{\bf
SL}_2\times{\bf SL}_2$ and its action on $\mathfrak
u_{l-2,1}^-$ is equivalent to that determined by the
pair $({\bf Sp}_{l-2}\times{\bf SL}_2\times{\bf SL}_2,
\Lambda_1\otimes\Lambda_1 \otimes\Lambda_1)$. As, by
\cite[\S\,6]{sk}, this pair is not prehomogeneous, we
conclude that if $l$ is even, then the action of
$L_{l-2}$ on $\mathfrak u_{l-2}^-$ is not locally
transitive.

Assume now that $l$ is odd. Then, according to
\cite[Proposition\,1.34\,(2)]{kkti}, the
prehomo\-ge\-nei\-ty of $({\bf GL}_{l-2}\times{\bf
SL}_2\times{\bf SL}_2,
\Lambda_1\otimes\Lambda_1\otimes\Lambda_1+
\Lambda_2\otimes 1\otimes 1)$ is equivalent to that~of
$({\bf GL}_{3}\times{\bf SL}_2\times{\bf SL}_2,
\Lambda_1\otimes\Lambda_1\otimes\Lambda_1+
\Lambda_2\otimes 1\otimes 1)$. The action determined
by $({\bf GL}_{3}\times{\bf SL}_2\times{\bf SL}_2,
\Lambda_1\otimes\Lambda_1\otimes\Lambda_1)$ is
equivalent to that of ${\bf GL}_3\times {\bf SO}_4$ on
${\rm Mat}_{3\times 4}$ given by
\begin{equation*}
g\cdot X:=AX{}^t\hskip -.7mm B,\quad g=(A, B)\in {\bf
GL}_3\times {\bf SO}_4,\ X\in {\rm Mat}_{3\times 4}.
\end{equation*}

By \cite[p.\,109]{sk}, the ${\bf GL}_3\times {\bf
SO}_4$-orbit of the matrix $[I_3\hskip .5mm 0]$ is
open in ${\rm Mat}_{3\times 4}$ and its stabilizer is
$$\bigl\{\bigl(A,
\Bigl[
\begin{matrix}{}^t\hskip -1.2mmA^{-1}&0\\
0&a\end{matrix} \Bigr] \bigr)\mid A\in {\bf O}_3,\
a=\pm 1,\ a\,{\rm det} A=1\bigr\}.$$ Hence the action
of the identity component of this stabilizer on the
second summand of the pair $({\bf GL}_{3}\times{\bf
SL}_2\times{\bf SL}_2,
\Lambda_1\otimes\Lambda_1\otimes\Lambda_1+
\Lambda_2\otimes 1\otimes 1)$ is equivalent to the
action of ${\bf SO}_3$ determined by the second
exterior power of its natural $3$-dimensional
representation. Since the last action is clearly not
locally transitive, Lemma~\ref{reduction} yields that
for odd $l$ the action of $L_{l-2}$ on $\mathfrak
u_{l-2}^-$ is not locally transitive.

Summing up, we obtain that for every $l$ the action of
$L_{l-2}$ on $\mathfrak u_{l-2}^-$ is not locally
transitive.

\smallskip

{\it Step} 6.\;Let $2<i<l-2$. From \eqref{diagrdi>1}
we deduce that the action of $L_{i}$ on $\mathfrak
u_{i}^-$ is equivalent to that determined by the pair
$({\bf GL}_i\times {\bf SO}_{2(l-i)},
\Lambda_1\otimes\Lambda_1+\Lambda_2\otimes 1)$. If the
action of $L_i$ on $\mathfrak u_i^-$ would be locally
transitive, then all the more the action determined by
the pair $({\bf GL}_1^2\times{\bf SL}_i\times {\bf
SO}_{2(l-i)},
\Lambda_1\otimes\Lambda_1+\Lambda_2\otimes 1)$, where
${\bf GL}_1^2$ acts on the summands by independent
scalar multiplications, would be locally transitive.
In turn, this would mean that if $2(l-i)>i$, then the
last pair is $2$-simple prehomogeneous of type I in
the sense of \cite{kkiy}. However the classification
of such pairs obtained in \cite[Section~3]{kkiy} shows
that it is not so. Hence, for $2(l-i)>i$, the action
of $L_i$ on $\mathfrak u_i^-$ is not locally
transitive.

Assume now that $2(l-i)\leqslant i$. Then, by the same
argument, if the action of $L_i$ on $\mathfrak u_i^-$
would be locally transitive, then the pair $({\bf
GL}_1^2\times{\bf SL}_i\times {\bf SO}_{2(l-i)},
\Lambda_1\otimes\Lambda_1+\Lambda_2\otimes 1)$ would
be $2$-simple prehomogeneous of type II in the sense
of \cite{kkiy}, \cite{kkti}. However the
classification of such pairs obtained in
\cite[Section\,5]{kkti} shows that it is not so. Hence
if $2(l-i)\leqslant i$, then the action of $L_i$ on
$\mathfrak u_i^-$ is not locally transitive.

Summing up, we obtain that if $2<i<l-2$, then the
action of $L_i$ on $\mathfrak u_i^-$ is not locally
transitive. \quad $\square$
\renewcommand{\qed}{}\end{proof}

\section{\bf \bf \boldmath ${\rm gtd}(L_i:\mathfrak  u_i^-)$
for \boldmath$G$ of type \boldmath${\sf
E}_6$}\label{E6}

\begin{proposition}\label{e6} Let $G$ be
a connected simple algebraic group of type ${\sf
E}_6$. Then
\begin{equation*}
{\rm gtd}(L_i:\mathfrak u_i^-)=\begin{cases}0& \mbox
{if}\hskip 2mm i=2, 4,\\ 1&\mbox{if}\hskip 2mm
i=3, 5,\\
2& \mbox{if}\hskip 2mm i=1, 6. \end{cases}
\end{equation*}
\end{proposition}
\begin{proof}
{\it Step} 1.\;Let $i=4$.
We have
$\dim\,G=78$. The type of $(L_4, L_4)$ is ${\sf
A}_1+{\sf A}_2+{\sf A}_2$, so $\dim\,L_4=20$,
$\dim\,\mathfrak u_4^-=29$. As
$\dim\,L_4<\dim\,\mathfrak u_4^-$, the action of $L_4$
on $\mathfrak u_4^-$ is not locally transitive.

\smallskip

{\it Step} 2.\;Let $i=1, 6$. By Remark \ref{abel} and
Corollary~\ref{abelianrad} of Proposition~\ref{p-},
the action of $L_1$ on $\mathfrak  u_1^-$ locally
transitive.  The type of $(L_1, L_1)$ is ${\sf D}_5$,
so $\dim\,L_1=46$, $\dim\,\mathfrak u_1^-=16$. As
$\mathfrak  u_1^-$ is abelian, the action of $L_1$ on
$\mathfrak  u_1^-$ is irreducible with the highest
weight $-\alpha_1=(-2, 0, 1, 0, 0, 0)$. Hence
$\mathfrak  u_1^-$ is a half-spinor module of $(L_1,
L_1)$. From \cite[Section 3, A, (17), (iii) and
Proposition 2.23]{kimu1}, \cite[Proposition 32]{sk} it
now follows that ${\rm gtd}(L_1:\mathfrak u_1^-)=2$.

By Proposition~\ref{p-} and \eqref{pp}, we have ${\rm
gtd}(L_6:\mathfrak u_6^-)=2$.

\smallskip

{\it Step} 3.\;Let $i=2$. The type of $(L_2, L_2)$ is
${\sf A}_5$, so $\dim L_2=36$, $\dim \mathfrak
u_2^-=21$. By inspection  of $\Phi^+$ in \cite[Planche
V]{bourbaki} we obtain  that, for the action of $L_2$
on $\mathfrak u_2^-$, there are exactly two shapes
$\alpha_2$ and $2\alpha_2$, and they determine the
highest roots $-\alpha_2=(0,-2,0,1,0,0)$ and
$-\alpha_1-2\alpha_2-2\alpha_3
-3\alpha_4-2\alpha_5-\alpha_6=(0,-1,0,0,0,0)$. Hence
$\mathfrak  u_2^-$ is the direct sum of two
irreducible $L_2$-modules $\mathfrak u_{21}^-$ and
$\mathfrak u_{22}^-$, where $\mathfrak u_{21}^-$,
considered as (irreducible) $(L_2, L_2)$-module,
 has the highest
weight with the Dynkin diagram
$$ \xymatrix@=4.5mm@M=
-.4mm{\circ\ar@{-}[r]
&\circ\ar@{-}[r]&\circ\ar@{-}[r]^{\hskip - 5.8mm 1}&
\circ\ar@{-}[r]& \circ}
$$ and $\mathfrak
u_{22}^-$ is a trivial one-dimensional  module. As the
action of $L_2$ on $\mathfrak u_{22}^-$ is locally
transitive, $(L_2)_z^0=(L_2, L_2)$ for a nonzero point
$z\in \mathfrak u_{22}^-$. Since the action of $(L_2,
L_2)$ on $\mathfrak u_{21}^-$ is not locally
transitive, \cite[\S7]{sk}, Lemma~\ref{reduction}
yields that the action of $L_2$ on $\mathfrak  u_2^-$
is not locally transitive.

\smallskip

{\it Step} 4.\;Let $i=3, 5$. The type of $(L_3, L_3)$
is ${\sf A}_1+{\sf A_4}$, so $\dim L_3= 28$, $\dim
\mathfrak u_3^-=25$. By inspection of $\Phi^+$ in
\cite[Planche V]{bourbaki} we obtain  that,  for the
action of $L_3$ on $\mathfrak u_3^-$, there are
exactly two shapes $\alpha_3$ and $2\alpha_3$, and
they determine the highest roots
$-\alpha_3=(1,0,-2,1,0,0)$ and
$-\alpha_1-\alpha_2-2\alpha_3-
2\alpha_4-\alpha_5=(0,0,-1,0,0,1)$. Hence $\mathfrak
u_3^-$ is the direct sum of two irreducible
$L_3$-modules $\mathfrak u_{31}^-$ and $\mathfrak
u_{32}^-$ that, considered as (irreducible) $(L_3,
L_3)$-modules, have respectively the highest weights
with the Dynkin diagrams
$$ \xymatrix@=4.5mm@M=
-.4mm{\circ&\circ\ar@{-}[r]^{\hskip - 17mm 1}&
\circ\ar@{-}[r]^{\hskip - 5.8mm 1}& \circ\ar@{-}[r]&
\circ}\hskip 6mm \text{and} \hskip 6mm
\xymatrix@=4.5mm@M= -.4mm{\circ&\circ\ar@{-}[r]&\circ
\ar@{-}[r]& \circ\ar@{-}[r]^{\hskip 5.7mm 1}&
\circ}\quad .$$

So we may (and shall) identify the Lie algebra of
$(L_3, L_3)$ with the Lie algebra of matrices
\begin{equation}\label{liel3l3}
\bigl\{
\Bigl[\begin{matrix}A&0\\
0&B
\end{matrix}\Bigr]\mid A\in {\rm
Mat}_{2\times 2},\ B\in {\rm Mat}_{5\times 5},\ {\rm
tr} A={\rm tr} B=0 \bigr\}
\end{equation}
(the Lie bracket is given by the commutator) and
$\mathfrak  u_{32}^-$ with the coordinate space $k^5$
on which this Lie algebra acts by the rule
\begin{equation*}
\begin{bmatrix}A&0\\0&B
\end{bmatrix}\cdot
v:=-{}^t\!Bv.
\end{equation*}
Then, by \cite[Lemma~1.4]{kkiy}, the open $L_3$-orbit
in $\mathfrak u_{31}^-$ contains a point $z$ such that
$$
\bigl\{\Bigl[\begin{matrix}A&0\\
0&B
\end{matrix}\Bigr]\mid A\!=\!
\Bigl[\begin{matrix}a_1&a_2\\
a_3&-a_1\end{matrix}\Bigr],
B\!=\!\Bigl[\begin{matrix}C&0\\D&A
\end{matrix}\Bigr],
C\!=\!-\!\begin{bmatrix}2a_1&2a_3&0\\
a_2&0&a_3\\
0&2a_2&-2a_1\end{bmatrix}\hskip -1mm,
D\!=\!\Bigl[\begin{matrix}a_4&a_5&a_6\\
a_5&a_6&a_7\end{matrix}\Bigr], a_i\!\in\! k\bigr\}
$$
is the Lie algebra of $(L_3, L_3)_z$. Hence $\dim
(L_3, L_3)_z=7$ and the Lie algebra of the $(L_3,
L_3)_z$-stabilizer of the point
$v:={}^t\!(1,0,0,0,0)\in \mathfrak u_{32}^-$ consists
of all $\bigl[\begin{smallmatrix}A&0\\0&B
\end{smallmatrix}\bigr]$ with
$a_1=a_3=0$. This shows that the $(L_3,L_3)_z$-orbit
of $v$ is $5$-dimensional and hence open in $\mathfrak
u_{32}^-$. All the more the action of $(L_3)_z$ on
$\mathfrak u_{32}^-$ is locally transitive.
Lemma~\ref{reduction} now yields that the action of
$L_3$ on $\mathfrak u_3^-$ is locally transitive. As
$\dim L_3<2\dim \mathfrak u_3^-$, this yields ${\rm
gtd}(L_3:\mathfrak u_3^-)=1$.

 By Proposition~\ref{p-} and \eqref{pp},
 we have ${\rm ltd}(L_5:\mathfrak
u_5^-)=1$. \quad $\square$
\renewcommand{\qed}{}\end{proof}

\section{\bf \bf \boldmath ${\rm gtd}(L_i:\mathfrak  u_i^-)$
for \boldmath$G$ of type \boldmath${\sf
E}_7$}\label{ge7}

\begin{proposition}\label{e7} Let $G$ be
a connected simple algebraic group of type ${\sf
E}_7$. Then
\begin{equation*}
{\rm gtd}(L_i:\mathfrak u_i^-)=\begin{cases}0& \mbox
{if}\hskip 2mm i\neq 7,\\ 1&\mbox{if}\hskip 2mm
i=7.\end{cases}
\end{equation*}
\end{proposition}

\begin{proof}
{\it Step} 1.\;We have $\dim\,G=133$. Let $i=3, 4, 5$.
Then the type of $(L_i, L_i)$ is respectively
 ${\sf A}_1+{\sf A}_5$, ${\sf
A}_1+{\sf A}_2+{\sf A}_3$, ${\sf A}_2+{\sf A}_4$.
Hence respectively $\dim\,L_i=39, 27, 33$ and
$\dim\,\mathfrak u_i^-=47, 53, 50$. As
$\dim\,L_i<\dim\,\mathfrak u_i^-$, the action of $L_i$
on $\mathfrak u_i^-$ is not locally transitive.

\smallskip

{\it Step}  2.\;Let $i=1$. The type of $(L_1, L_1)$ is
${\sf D}_6$, so $\dim\,L_1=67$, $\dim\,\mathfrak
u_1^-=33$. By inspection  of $\Phi^+$ in \cite[Planche
VI]{bourbaki} we obtain that,  for the action of $L_1$
on $\mathfrak u_1^-$, there are exactly two shapes
$\alpha_1$ and $2\alpha_1$, and they determine the
highest weights $-\alpha_1=(-2,0,1,0,0,0,0)$ and
$-2\alpha_1-2\alpha_2-3\alpha_3-4\alpha_4-
3\alpha_5-2\alpha_6-\alpha_7=(-1, 0, 0, 0, 0, 0, 0)$.
Hence $\mathfrak  u_1^-$ is the direct sum  of two
irreducible $L_1$-modules $\mathfrak  u_{11}^-$ and
$\mathfrak  u_{12}^-$ that, considered as
(irreducible) $(L_1, L_1)$-modules,
 are respectively a
half-spinor and a trivial $1$-dimen\-sio\-nal module.
As the action of $L_1$ on $\mathfrak u_{12}^-$ is
locally transitive, $(L_1)_z^0=(L_1, L_1)$ for a
general point $z\in \mathfrak  u_{12}^-$. Since the
action of $(L_1, L_1)$ on $\mathfrak u_{11}^-$ is not
locally transitive, \cite[\S\,7]{sk},
Lemma~\ref{reduction} yields that the action of $L_1$
on $\mathfrak  u_1^-$ is not locally transitive.

\smallskip

{\it Step}  3.\;Let $i=2$. The type of $(L_2, L_2)$ is
${\sf A}_6$, so $\dim\,L_2=49$, $\dim\,\mathfrak
u_2^-=42$. By inspection  of $\Phi^+$ in \cite[Planche
VI]{bourbaki} we obtain that,  for the action of $L_2$
on $\mathfrak u_2^-$, there are exactly two shapes
$\alpha_2$ and $2\alpha_2$, and they determine the
highest weights $-\alpha_2=(0, -2, 0, 1, 0, 0, 0)$ and
$-\alpha_1-2\alpha_2-2\alpha_3
-3\alpha_4-2\alpha_5-\alpha_6= (0,-1,0,0,0,0,1)$.
Hence $\mathfrak u_2^-$ is the direct sum  of two
irreducible $L_2$-modules $\mathfrak u_{21}^-$ and
$\mathfrak  u_{22}^-$ that, considered as
(irreducible) $(L_2, L_2)$-modules, have respectively
the highest weights with the Dynkin diagrams
$$ \xymatrix@=4.5mm@M=
-.4mm{\circ\ar@{-}[r]
&\circ\ar@{-}[r]&\circ\ar@{-}[r]^{\hskip - 5.8mm 1}&
\circ\ar@{-}[r]& \circ\ar@{-}[r]&\circ}\hskip 6mm
\text{and} \hskip 6mm \xymatrix@=4.5mm@M=
-.4mm{\circ\ar@{-}[r]
&\circ\ar@{-}[r]&\circ\ar@{-}[r]& \circ\ar@{-}[r]&
\circ\ar@{-}[r]^{\hskip 5.6mm 1}&\circ}\quad .$$

By \cite[\S\,7, I, (6)]{sk}, if $z$ is a general point
of $\mathfrak u_{21}^-$, then $(L_2)_z^0$ is a simple
algebraic group of type ${\sf G}_2$. Hence
$(L_2)_z^0\subset (L_2, L_2)$ and, as $(L_2, L_2)$ is
simple, the action of $(L_2)_z^0$ on $\mathfrak
u_{22}^-$ is nontrivial. Since the dimension of every
nontrivial module of a simple group of type ${\sf
G}_2$ is at least $7=\dim \mathfrak u_{22}^-$, and, by
\cite[\S\,7]{sk}, the action of ${\bf G}_2$ on every
such module is not locally transitive, this yields
that the action of $(L_2)_z^0$ on $\mathfrak u_{22}^-$
is
 not locally transitive. From
Lemma~\ref{reduction} we then deduce that the action
of $L_2$ on $\mathfrak u_2^-$ is not locally
transitive as well.

\smallskip

{\it Step} 4.\;Let $i=6$. The type of $(L_6, L_6)$ is
${\sf A}_1+{\sf D}_5$, so $\dim L_6=49$, $\dim
\mathfrak u_6^-=42$. By inspection of $\Phi^+$ in
\cite[Planche VI]{bourbaki} we obtain  that,  for the
action of $L_6$ on $\mathfrak u_6^-$, there are
exactly two shapes $\alpha_6$ and $2\alpha_6$, and
they determine the highest weights $-\alpha_6= (0, 0,
0, 0, 1, -2, 1)$ and
$-\alpha_2-\alpha_3-2\alpha_4-2\alpha_5-2\alpha_6
-\alpha_7=(1, 0, 0, 0, 0, -1, 0)$. Hence $\mathfrak
u_6^-$ is the direct sum of two irreducible
$L_2$-modules $\mathfrak u_{61}^-$ and $\mathfrak
u_{62}^-$ that, considered as (irreducible) $(L_2,
L_2)$-modules, have respectively the highest weights
with the Dynkin diagrams
$$
\hskip -20mm
\begin{matrix}
\xymatrix@=4.5mm@M=
-.4mm@R=3mm{ &&&\circ&\\
\circ\ar@{-}[r] &\circ\ar@{-}[r]^{\hskip 28mm 1}&
\circ\ar@{-}[ur]^{\hskip 38mm 1}\ar@{-}[dr]&&\circ\\
&&&\circ& }
\end{matrix}
\hskip 10mm \text{and} \hskip 10mm
\begin{matrix}
\xymatrix@=4.5mm@M=
-.4mm@R=3mm{ &&&\circ&\\
\circ\ar@{-}[r]^{\hskip -5.8mm 1} &\circ\ar@{-}[r]&
\circ\ar@{-}[ur]\ar@{-}[dr]&&\circ\\
&&&\circ& }
\end{matrix}\quad .
$$

If the action of $L_6$ on $\mathfrak u_6^-$ would be
locally transitive, then all the more the action of
${\bf GL}_1\times L_6$ on $\mathfrak u_6^-$, where the
first factor acts on $\mathfrak  u_{61}^-$ by scalar
multiplication and trivially on $\mathfrak u_{62}^-$,
would be locally transitive. Using the notation and
terminology of \cite{kkiy}, this in turn would mean
that the pair $({\bf GL}_1^2\times {\bf
Spin}_{10}\times {\bf SL}_2,
\Lambda'\otimes\Lambda_1+\Lambda_1\otimes 1)$  is
$2$-simple prehomogeneous of type I. However the
classification of such pairs obtained in \cite[Section
3]{kkiy} shows that it is not so. Thus the action of
$L_6$ on $\mathfrak u_6^-$ is not locally transitive.

\smallskip

{\it Step} 5.\;Let $i=7$. By Remark \ref{abel} and
Corollary~\ref{abelianrad} of Proposition~\ref{p-},
the action of $L_7$ on $\mathfrak  u_7^-$ is locally
transitive. The type of $(L_7, L_7)$ is ${\sf E}_6$,
so $\dim\,L_7=79$, $\dim\,\mathfrak  u_7^-=27$. By the
dimension reason, $\mathfrak  u_7^-$ is a minimal
irreducible $(L_7, L_7)$-module. Hence \cite[\S\,7, I,
(27)]{sk} yields that $(L_7)_z^0$ for a general point
$z\in \mathfrak u_7^-$ is a simple algebraic group of
type ${\sf F}_4$. Since the dimension of every
nontrivial module of a simple group of type ${\sf
F}_4$ is at least $26$,  the $(L_7)_z^0$-module
$\mathfrak  u_7^-$ contains a trivial one-dimensional
submodule and hence is not locally transitive.
Lemma~\ref{reduction} now yields
 ${\rm ltd}(L_7:\mathfrak u_7^-)=1$. \quad $\square$
\renewcommand{\qed}{}\end{proof}

\section{\bf \bf \boldmath ${\rm gtd}(L_i:\mathfrak  u_i^-)$
for \boldmath$G$ of type \boldmath${\sf
E}_8$}\label{E8}

\begin{proposition}\label{e8} Let $G$ be
a connected simple algebraic group of type ${\sf
E}_8$. Then
\begin{equation*}
{\rm gtd}(L_i:\mathfrak u_i^-)=0\ \mbox{\ for every \
} i.
\end{equation*}
\end{proposition}

\begin{proof}
{\it Step} 1.\;We have $\dim\,G=248$. Let $i=2, 3, 4,
5, 6, 7$. Then the type of $(L_i, L_i)$ is
respecti\-ve\-ly ${\sf A}_7$, ${\sc A}_1+{\sf A}_6$,
${\sf A}_1+{\sf A}_2+ {\sf A}_4$, ${\sf A}_3+{\sf
A}_4$, ${\sf A}_2+{\sf D}_5$, ${\sf A}_1+{\sf E}_6$.
Hence
 respectively $\dim\,L_i=64, 52, 36, 40,
54, 82$ and $\dim\,\mathfrak u_i^-=92, 98, 106, 104,
97, 83$. As $\dim\,L_i<\dim\,\mathfrak u_i^-$, the
action of $L_i$ on $\mathfrak u_i^-$ is not locally
transitive.

\smallskip

{\it Step} 2.\;Let $i=1$. The type of $(L_1, L_1)$ is
${\sf D}_7$, so $\dim L_1=92$, $\dim \mathfrak
u_1^-=78$. By inspection  of $\Phi^+$ in \cite[Planche
VII]{bourbaki}  we obtain that,  for the action of
$L_1$ on $\mathfrak u_1^-$, there are exactly two
shapes $\alpha_6$ and $2\alpha_6$, and they determine
the highest weights $-\alpha_1=(-2, 0, 1, 0, 0, 0, 0,
0)$ and $-\alpha_1-2\alpha_2
-3\alpha_3-4\alpha_4-3\alpha_5-2\alpha_6
-\alpha_7=(-1, 0, 0, 0, 0, 0, 0, 1)$. Hence $\mathfrak
u_1^-$ is the direct sum of two irreducible
$L_1$-modules $\mathfrak u_{11}^-$ and $\mathfrak
u_{12}^-$ that, considered as (irreducible) $(L_1,
L_1)$-modules, have respectively the highest weights
with the Dynkin diagrams
$$
\hskip 5mm
\begin{matrix}
\xymatrix@=4.5mm@M= -.4mm@R=3mm{
\circ\ar@{-}[dr]^{\hskip - 7mm
1}&&&&&\\
&\circ\ar@{-}[r]& \circ\ar@{-}[r]&
\circ\ar@{-}[r]&\circ\ar@{-}[r]&\circ\\
\circ\ar@{-}[ur]&&&&& }
\end{matrix}
\hskip 10mm \text{and} \hskip 10mm
\begin{matrix}
\xymatrix@=4.5mm@M=
-.4mm@R=3mm{ \circ&&&&&\\
&\circ\ar@{-}[ul]\ar@{-}[dl]\ar@{-}[r]&
\circ\ar@{-}[r]&
\circ\ar@{-}[r]&\circ\ar@{-}[r]^{\hskip
5.5mm 1}&\circ\\
\circ&&&&& }
\end{matrix}
\quad .
$$

If the action of $L_1$ on $\mathfrak u_1^-$ would be
locally transitive, then all the more the action of
${\bf G}_m\times L_1$ on $\mathfrak u_1^-$, where the
first factor acts on $\mathfrak  u_{11}^-$ by scalar
multiplication and trivially on $\mathfrak u_{12}^-$,
would be locally transitive. This in turn would mean
that, in the notation and terminology of \cite{kimu1},
$({\bf GL}_1^2\times {\bf Spin}_{14},
\Lambda'+\Lambda_1)$ is a prehomogeneous vector space
with scalar multiplications. However the
classification of such spaces obtained in
\cite[Section~3]{kimu1} shows that it is not so. Thus
the action of $L_6$ on $\mathfrak  u_6^-$ is not
locally transitive.

\smallskip

{\it Step} 3.\;Let $i=8$. The type of $(L_8, L_8)$ is
${\sf E}_7$, so $\dim L_8=134$, $\dim u_8^-=57$. By
inspection  of $\Phi^+$ in \cite[Planche
VII]{bourbaki}  we obtain that,  for the action of
$L_8$ on $\mathfrak u_8^-$, there are exactly two
shapes $\alpha_8$ and $2\alpha_8$, and they determine
the highest weights $-\alpha_8=(0, 0, 0, 0, 0, 0, 1,
-2)$ and $-2\alpha_8=(0, 0, 0, 0, 0, 0, 0, -1)$. Hence
$\mathfrak u_8^-$ is the direct sum  of two
irreducible $L_8$-modules $\mathfrak  u_{81}^-$ and
$\mathfrak u_{82}^-$ that, considered as (irreducible)
$(L_8, L_8)$-modules,
 are
 respectively the unique
 $56$-dimensional
and a trivial $1$-dimen\-sio\-nal module. As the
action of $(L_8, L_8)$ on $\mathfrak u_{81}^-$ is not
locally transitive, \cite[\S\,7]{sk}, the same
argument as in Subsection~\ref{ge7} for $L_1$ shows
that the action of $L_8$ on $\mathfrak  u_8^-$ is not
locally transitive. \quad $\square$
\renewcommand{\qed}{}\end{proof}

\section{\bf \bf \boldmath ${\rm gtd}(L_i:\mathfrak  u_i^-)$
for \boldmath$G$ of type \boldmath${\sf
F}_4$}\label{F4}

\begin{proposition}\label{f4} Let $G$ be
a connected simple algebraic group of type ${\sf
F}_4$. Then
\begin{equation*}
{\rm gtd}(L_i:\mathfrak u_i^-)=0\ \mbox{\ for every \
} i.
\end{equation*}
\end{proposition}

\begin{proof}
{\it Step} 1.\;We have $\dim\,G=52$. For $i=2, 3$ the
type of $(L_i, L_i)$ is ${\sf A}_1+{\sf A_2}$, so
$\dim L_i=12$, $\dim \mathfrak u_i^-=20$. As
$\dim\,L_i<\dim\,\mathfrak u_i^-$, the action of $L_i$
on $\mathfrak u_i^-$ is not locally transitive.

\smallskip

{\it Step} 2.\;Let $i=1$. The type of $(L_1, L_1)$ is
${\sf C}_3$, so $\dim L_1=22$ and $\dim \mathfrak
u_1^-=15$. By inspection of $\Phi^+$ in \cite[Planche
VIII]{bourbaki}  we obtain  that,  for the action of
$L_1$ on $\mathfrak u_1^-$, there are exactly two
shapes $\alpha_1$ and $2\alpha_1$, and they determine
the highest weights $-\alpha_1=(-2, 1, 0, 0)$ and
$-2\alpha_1-3\alpha_2-4\alpha_3-2\alpha_4= (-1, 0, 0,
0)$. Hence $\mathfrak u_1^-$ is the direct sum  of two
irreducible $L_1$-modules $\mathfrak  u_{11}^-$ and
$\mathfrak  u_{12}^-$, where $\mathfrak u_{11}^-$,
considered as (irreducible) $(L_1, L_1)$-module,  has
the highest weight with the Dynkin diagram
$$ \xymatrix@=7mm
@M= -.4mm{\circ\ar@{=>}[r]^{\hskip - 8.1mm 1}
&\circ\ar@{-}[r]&\circ}$$ \noindent and $\mathfrak
u_{12}^-$ is a trivial $1$-dimen\-sional module. As
the action of $(L_1, L_1)$ on $\mathfrak u_{11}^-$ is
not locally transitive, \cite[\S\,7]{sk}, the same
argument as in Subsection~\ref{ge7} for $L_1$ shows
that the action of $L_1$ on $\mathfrak  u_1^-$ is not
locally transitive.

\smallskip

{\it Step} 3.\;Let $i=4$. The type of $(L_4, L_4)$ is
${\sf B}_3$, so $\dim L_1=22$ and $\dim \mathfrak
u_1^-=15$. By inspection of $\Phi^+$ in \cite[Planche
VIII]{bourbaki}  we obtain  that,  for the action of
$L_4$ on $\mathfrak u_4^-$, there are exactly two
shapes $\alpha_4$ and $2\alpha_4$, and they determine
the highest weights $-\alpha_4=(0, 0, 1, -2)$ and
$-\alpha_2-2\alpha_3-2\alpha_4=(1, 0, 0, -2)$. Hence
$\mathfrak  u_4^-$ is the direct sum  of two
irreducible $L_4$-modules $\mathfrak  u_{41}^-$ and
$\mathfrak u_{42}^-$ that, considered as (irreducible)
$(L_4, L_4)$-modules, have respectively the highest
weights with the Dynkin diagrams
$$
\begin{matrix}
\hskip 5mm \xymatrix@=7mm @M= -.4mm{\circ\ar@{-}[r]
&\circ\ar@{=>}[r]^{\hskip 8.1mm 1} &\circ}
\end{matrix}
\hskip 10mm \text{and} \hskip 10mm \quad
\begin{matrix}
\xymatrix@=7mm @M= -.4mm{\circ\ar@{-}[r]^{\hskip -
8.1mm 1} &\circ\ar@{=>}[r]&\circ}
\end{matrix}
\quad .
$$

 From \cite[\S\,7, I,
(16)]{sk} we now deduce that $(L_4)_z^0$ for a general
point $z\in \mathfrak u_{41}^-$ is a simple algebraic
group of type ${\sf G}_2$. Hence $(L_4)_z^0\subset
(L_4, L_4)$. As the action of $(L_4, L_4)$ on
$\mathfrak  u_{42}^-$ is clearly not locally
transitive, \cite[\S\,7]{sk}, this yields that the
action of $(L_4)_z^0$ on $\mathfrak  u_{42}^-$ is not
locally transitive as well. From Lemma~\ref{reduction}
we then deduce that the action of $L_4$ on $\mathfrak
u_4^-$ is not locally transitive. \quad $\square$
\renewcommand{\qed}{}\end{proof}

\section{\bf \bf \boldmath ${\rm gtd}(L_i:\mathfrak  u_i^-)$
for \boldmath$G$ of type \boldmath${\sf
G}_2$}\label{G2}

\begin{proposition}\label{g2} Let $G$ be
a connected simple algebraic group of type ${\sf
G}_2$. Then
\begin{equation*}
{\rm gtd}(L_i:\mathfrak u_i^-)=0\ \mbox{\ for every \
} i.
\end{equation*}
\end{proposition}
\begin{proof}
We have $\dim\,G=14$. For every $i$ the type of $(L_i,
L_i)$ is ${\sf A}_1$, so $\dim L_i=4$, $\dim \mathfrak
u_i^-=5$. As $\dim\,L_i<\dim\,\mathfrak u_i^-$, the
action of $L_i$ on $\mathfrak u_i^-$ is not locally
transitive. \quad $\square$
\renewcommand{\qed}{}\end{proof}

\section{\bf Proofs of
Theorems~\ref{main1}--\ref{main6}}

\noindent{\it Proof of Theorem~{\rm \ref{main5}.}}
Statement (i) follows from Proposition~\ref{p-}(iii),
and  (ii) from Propositions~\ref{al}--\ref{g2}. \quad
$\square$

\bigskip

\noindent{\it Proof of Theorem~{\rm \ref{main2}.}} The
claim follows from Proposition~\ref{1-2}(ii). \quad
$\square$

\bigskip

\noindent{\it Proof of Theorem~{\rm \ref{main4}.}}
Statement (a) follows from Proposition~\ref{p-}, and
 (b) from Theorem~\ref{main5}. \quad
$\square$

\bigskip

\noindent{\it Proof of Theorem~{\rm \ref{main6}.}} The
argument is based on the following facts proved in
\cite{popov-vinberg1}. Let $\mathcal O(\varpi)$ be the
$G$-orbit of a nonzero $B$-semi-invariant vector in
$E(\varpi)$, and let $\mathcal C(\varpi)$ be the
closure of $\mathcal O(\varpi)$ in $E(\varpi)$. Then
$\mathcal C(\varpi)$ is a cone, i.e., stable with
respect to the action of ${\bf G}_m$ on $E(\varpi)$ by
scalar multiplications, and $\mathcal
O(\varpi)=\mathcal C(\varpi)\setminus \{0\}$,
\cite[Theorem~1]{popov-vinberg1}.
 This ${\bf G}_m$-action
commutes with the $G$-action and yields a $G$-stable
$\mathbb Z_+$-grading of the algebra $k[\mathcal
C(\varpi)]$,
\begin{equation}\label{grad}\textstyle
k[\mathcal C(\varpi)]=\bigoplus_{n\in \mathbb
Z_+}k[\mathcal C(\varpi)]_n.\end{equation} For every
$n\in \mathbb Z_+$, the $G$-module $k[\mathcal
C(\varpi)]_n$ is isomorphic to $E(n\varpi^*)$,
\cite[Theorem 2]{popov-vinberg1}. If $\varpi$ is
dominant, then  $k[\mathcal C(\varpi)]$ is a unique
factorization domain, \cite[Theorem
4]{popov-vinberg1}.

Thus $G\times {\bf G}_m^d$ acts on $\mathcal
C(\varpi)^d$, and $\mathcal O(\varpi)^d$ is an open
$G\times {\bf G}_m^d$-stable subset of $\mathcal
C(\varpi)^d$. Restricting the action to ${\bf G}_m^d$
yields a $G$-stable $\mathbb Z_+^d$-grading of the
algebra $k[\mathcal C(\varpi)^d]$,
\begin{equation}\label{mgrad}
\textstyle k[\mathcal C(\varpi)^d]=
\bigoplus_{(n_1,\ldots, n_d)\in \mathbb Z_+^d}
k[\mathcal C(\varpi)^d]_{(n_1,\ldots,
n_d)}.\end{equation} Since $k[\mathcal C(\varpi)^d]$
and $k[\mathcal C(\varpi)]^{\otimes d}$ are
isomorphic,
 \eqref{grad} and \eqref{mgrad} yield
 that the
$G$-modules $k[\mathcal C(\varpi)^d]_{(n_1,\ldots,
n_d)}$ and $E(n_1\varpi^*)\otimes\ldots\otimes
E(n_d\varpi^*)$ are isomorphic for every $(n_1,\ldots,
n_d)\in \mathbb Z_+^d$.

Let now $\pi_\varpi: \mathcal O(\varpi)\longrightarrow
G/P(\varpi)$ be the natural projection. The
$G$-equivariant morphism
\begin{equation*}
\pi_\varpi^d:\ \mathcal O(\varpi)^d\longrightarrow
\bigl(G/P(\varpi)\bigr)^d
\end{equation*}
is the quotient by  ${\bf G}_m^d$-action. Hence it
yields an isomorphism of invariant fields
\begin{equation}\label{fields1}
 k\bigl(\bigl(G/P(\varpi)\bigr)^d\bigr)^G
 \overset{\simeq}{\longrightarrow}
 k\bigl(\mathcal O(\varpi)^d\bigr)^{G\times {\bf
 G}_m^d}.
\end{equation}
Note that Definition~\ref{defgtdeg} and Rosenlicht's
theorem, \cite{rosenlicht}, yield the equivalence
\begin{equation}\label{fields2}
k\bigl(\bigl(G/P(\varpi)\bigr)^d\bigr)^G=k \
\Longleftrightarrow\ {\rm gtd} (G:
G/P(\varpi))\geqslant d.
\end{equation}
From \eqref{fields1} and \eqref{fields2} we deduce the
equivalence
\begin{equation}\label{fields3}
 k\bigl(\mathcal O(\varpi)^d\bigr)^{G\times {\bf
 G}_m^d}=k
 \
\Longleftrightarrow\ {\rm gtd} (G:
G/P(\varpi))\geqslant d.
\end{equation}
As $\mathcal O(\varpi)^d$ is open in $\mathcal
C(\varpi)^d$, we have $k(\mathcal
O(\varpi)^d)=k(\mathcal C(\varpi)^d)$. Hence
\eqref{fields3} yields
\begin{equation}\label{fields4}
 k\bigl(\mathcal C(\varpi)^d\bigr)^{G\times {\bf
 G}_m^d}=k
 \
\Longleftrightarrow\ {\rm gtd} (G:
G/P(\varpi))\geqslant d.
\end{equation}

We can now prove statements  (i) and (ii) of
Theorem~\ref{main6}.

(i)  Assume the contrary. Take $n_1,\ldots,n_d\in
\mathbb Z_+$ such that
\begin{equation}\label{tens}
\dim \bigl(E(n_1\varpi^*)\otimes\ldots \otimes
E(n_d\varpi^*)\bigr)^G\geqslant 2.
\end{equation}
Then $\dim k[\mathcal
C(\varpi)^d]_{(n_1,\ldots,n_d)}^G \geqslant 2$, so
there are nonzero elements $f_1, f_2\in k[\mathcal
C(\varpi)^d]_{(n_1,\ldots,n_d)}^G$ such that
$f_1/f_2\notin k$. Since $f_1$ and $f_2$ are ${\bf
G}_m^d$-semi-invariants of the same weight
$(n_1,\ldots n_d)$, we have $f_1/f_2\in
k\bigl(\mathcal C(\varpi)^d\bigr)^{G\times {\bf
 G}_m^d}$. Thus $k\bigl(\mathcal
 C(\varpi)^d\bigr)^{G\times {\bf
 G}_m^d}\neq k$. By \eqref{fields4}, this contradicts
the condition ${\rm gtd} (G: G/P(\varpi))\geqslant d$.

 \smallskip

(ii) Assume the contrary. Then by \eqref{fields4},
there is a nonconstant rational function $f\in
k\bigl(\mathcal
 C(\varpi)^d\bigr)^{G\times {\bf
 G}_m^d}$. Take now into account that (a) $k\bigl(\mathcal
 C(\varpi)^d\bigr)$ is the field of quotients of
$k[\mathcal C(\varpi)^d]$ (as $\mathcal
 C(\varpi)^d$ is affine); (b) $k[\mathcal
C(\varpi)^d]$ is a unique factorization domain (as
$\varpi$ is do\-mi\-nant); (c) $G\times {\bf G}_m^d$
is connected. By \cite[Theorem 3.3]{popov-vinberg2}
these properties yield that $f=f_1/f_2$ for some $f_1,
f_2\in k[\mathcal C(\varpi)^d]$ which are $G\times
{\bf G}_m^d$-semi-invariants of the same weight. Since
$G$ has no nontrivial characters, the latter means
that $f_1, f_2\in k[\mathcal
C(\varpi)^d]_{(n_1,\ldots,n_d)}^G$ for some
$n_1,\ldots, n_d\in \mathbb Z_+$. As $f_1/f_2\notin
k$, this yields $\dim k[\mathcal
C(\varpi)^d]_{(n_1,\ldots,n_d)}^G\geqslant 2$. Hence
\eqref{tens} holds, and this contradicts \eqref{<2}.
\quad $\square$

\bigskip

\noindent{\it Proof of Theorem~{\rm \ref{main3}.}}
{\it Step} 1.\;If $P_i$ is conjugate to $P_i^-$ (i.e.,
$\varepsilon(\alpha_i)=\alpha_i$, see
Section~\ref{stpar}), then the claim follows from
Proposition~\ref{p-}, its Corollary~\ref{cor}, and
Theorem~\ref{main5}. This covers all but the following
cases:
\begin{enumerate} \item[(a)] $G$ is of type ${\sf
A}_l$, and $2i\neq l+1$; \item[(b)] $G$ is of type
${\sf D}_l$, $l$ is odd, and $i=l-1, l$; \item[(c)]
$G$ is of type ${\sf E}_6$, and $i=1, 3, 5,
6$.\end{enumerate}

\medskip

{\it Step} 2.\;Consider case (b). By \eqref{pp}, we
have ${\rm gtd}(G: G/P_{l-1})={\rm gtd}(G: G/P_{l})$.
By Proposition~\ref{p-}, its Corollary~\ref{cor}, and
Theorem~\ref{main5},  we have ${\rm
gtd}(G:G/P_{l-1})\geqslant 3$.

Thus it suffices to prove ${\rm gtd}(G:G/P_{l-1})<4$.
Towards this end we apply Theorem~\ref{main6}. First,
note that for any semisimple $G$ and $\lambda, \mu\in
{\rm P}_{++}$, we have $E(\lambda)\otimes E(\mu)={\rm
Hom}(E(\lambda)^*, E(\mu))={\rm Hom}(E(\lambda^*),
E(\mu))$, and the elements of ${\rm Hom}(E(\lambda^*),
E(\mu))^G$ are precisely
 $G$-module homomorphisms
$E(\lambda^*)\rightarrow E(\mu)$. Since
$E(\lambda^*)$ and $E(\mu)$ are simple, this yields
\begin{equation}\label{multiplicity}
\dim\bigl(E(\lambda)\otimes
E(\mu)\bigr)^G=\begin{cases} 1&\mbox{if
$\lambda^*=\mu$},\\
0&\mbox{otherwise}.
\end{cases}
\end{equation}
In case (b), we have $\varpi^*_s=\varpi_s$ for every
$s\leqslant l-2$, whence by \eqref{multiplicity}
\begin{equation}\label{tensor1}
\dim \bigl(E(\varpi_s)\otimes E(\varpi_s)\bigr)^G=
1\quad\mbox{for every $s\leqslant l-2$}.
\end{equation}
On the other hand,  by \cite[Table 5]{oni-vinb}, we
have
\begin{gather}
\label{tensor2}
\textstyle E(\varpi_l)\otimes
E(\varpi_l)=E(2\varpi_l)\oplus\bigoplus_{i=1}^{\infty}
E(\varpi_{l-2i})
\end{gather}
where, by definition, $\varpi_t=0$ for $t<0$. Since
$l\geqslant 4$, from \eqref{tensor1} and
\eqref{tensor2} it then clearly follows that
\begin{equation} \label{tensor3}
\dim
\bigl(E(\varpi_l)\otimes E(\varpi_l)\otimes
E(\varpi_l)\otimes E(\varpi_l)\bigr)^G\geqslant 2.
\end{equation} Since
$E(\varpi_{l-1})^*=E(\varpi_l)$,
Theorem~\ref{main6}(a) and  \eqref{tensor3} now yield
${\rm gtd}(G:G/P_{l-1})<4$. This com\-p\-le\-tes the
proof in case (b).

\medskip

{\it Step} 3.\; Consider case (c). By \eqref{pp}, we
have ${\rm gtd}(G: G/P_{1})={\rm gtd}(G: G/P_{6})$ and
${\rm gtd}(G: G/P_{3})={\rm gtd}(G: G/P_{5})$.
Since $\dim G=78$, $\dim G/P_6=16$,
we have
$\dim (G/P_6)^5>\dim G$.
Hence ${\rm gtd}(G:G/P_6)\leqslant 4$.
By the Corollary of Lemma~\ref{parabolic}, we have
${\rm gtd}(G:G/P_5)\geqslant 2$.
 Thus it suffices to
prove that ${\rm gtd}(G:G/P_6)\geqslant 4$ and ${\rm
gtd}(G:G/P_5)< 3$.
Towards this end we apply Theorem~\ref{main6}.

Namely, we have
\begin{equation}\label{*}
\varpi_1^*=\varpi_6, \hskip 2mm \varpi_3^*=\varpi_5.
\end{equation}
By Theorem~\ref{main6} and \eqref{*}, proving ${\rm
gtd}(G:G/P_6)\geqslant 4$ is equivalent to proving
\begin{equation}\label{tensor4}
\dim \bigl(E(n_1\varpi_1)\otimes E(n_2\varpi_1)\otimes
E(n_3\varpi_1)\otimes E(n_4\varpi_1)\bigr)^G\leqslant
1\quad\mbox{for every $n_1,\ldots, n_4\in\mathbb
Z_+$}.
\end{equation}
To prove \eqref{tensor4}, we use that for every $r,
s\in\mathbb Z_+$ the following decomposition holds
 (see\break
\cite[1.3]{littelmann}):
\begin{equation}\label{tensor5}
E(r\varpi_1)\otimes E(s\varpi_1)=\hskip
-7mm\bigoplus_{\fontsize{8pt}{5mm}
\selectfont\biggl\{\begin{gathered}
a_1,\ldots, a_4\in\mathbb Z_+,\\[-8pt]
a_1+a_3+a_4=r,\\[-8pt]
a_2+a_3+a_4=s\end{gathered}
}\hskip -7mm
E\bigl((a_1+a_2)\varpi_1+a_3\varpi_3+a_4\varpi_6\bigr).
\end{equation}

\vskip -1mm \noindent Since, by \eqref{*}, we have
$\bigl((a_1+a_2)\varpi_1+a_3\varpi_3+a_4\varpi_6\bigr)^*=
a_4\varpi_1+a_3\varpi_5+(a_1+a_2)\varpi_6$, it follows
from \eqref{tensor5} and \eqref{multiplicity} that
$\dim \bigl(E(n_1\varpi_1)\otimes
E(n_2\varpi_1)\otimes E(n_3\varpi_1)\otimes
E(n_4\varpi_1)\bigr)^G$ is equal to the number of
solutions in $\mathbb Z_+$ of the following system of
eight linear equations in eight variables $a_1,\ldots,
a_4, b_1,\ldots, b_4$: \vskip -5mm
\begin{gather}
\left\{\hskip -59mm
\begin{align*}
a_4&=b_1+b_2,\\
a_3&=0,\\
b_3&=0,\\
a_1&+a_2=b_4,\\
a_1&+a_3+a_4=n_1,\\
a_2&+a_3+a_4=n_2,\\
b_1&+b_3+b_4=n_3,\\
b_2&+b_3+b_4=n_4.
\end{align*}
\right.
\end{gather}
Since this system is nondegenerate, there is at most
one such solution. Thus, \eqref{tensor4} holds; whence
${\rm gtd}(G:G/P_6)=4$.

By Theorem~\ref{main6} and \eqref{*}, proving ${\rm
gtd}(G: G/P_5)<3$ is equivalent to proving
\begin{equation}\label{tensor6}
\dim \bigl(E(n_1\varpi_3)\otimes E(n_2\varpi_3)\otimes
E(n_3\varpi_3)\bigr)^G\geqslant 2\quad\mbox{for some
$n_1, n_2, n_3\in\mathbb Z_+$}.
\end{equation}
Using Klimyk's formula, one checks  that the
decomposition of $E(2\varpi_3)\otimes E(2\varpi_3)$
into simple factors contains $E(2\varpi_5)$ with
multiplicity $2$ (using computer algebra system {\sf
LiE}, one obtains this decomposition in less than 1
second; this system is now available online  at {\tt
http:/\hskip -1.3mm/\break
wwwmathlabo.univ-poitiers.fr/$\sim$maavl/LiE/}).
Hence, by \eqref{*} and \eqref{multiplicity}, the
inequality \eqref{tensor6} holds for $n_1=n_2=n_3=2$.
This completes the proof in case (c).

\medskip

{\it Step} 4.\;To consider case (a), let now $G$ be of
type ${\sf A}_l$ and, more generally, no restrictions
are imposed on $i$. By Proposition~\ref{isogen}, we
may (and shall) assume that \eqref{slp} holds. Then
$G/P_i$ is the Grassmannian variety of $i$-dimensional
linear subspaces of $k^{l+1}$.

Given an integer $d\geqslant 1$, let ${\mathcal V}_d$
be the quiver with $d+1$ vertices, $d$ outside, one
inside, and the arrows from each vertex outside to a
vertex inside (the vertices are enumerated by
$1,\ldots, d+1$ so that the inside vertex is
enumerated by $1$):

\vskip 4mm

\begin{center}
\leavevmode \epsfxsize =3cm \epsffile{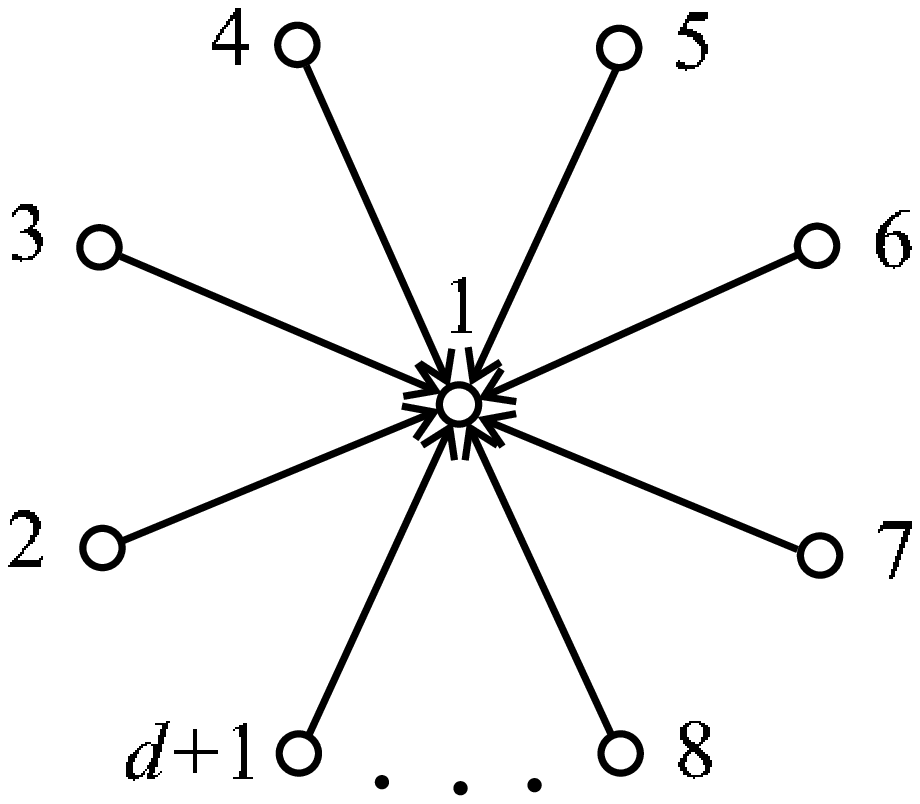}
\end{center}

\vskip 3mm

\noindent
 Given a vector
\begin{equation*}
\alpha:=(a_1,\ldots, a_{d+1})\in \mathbb Z_+^{d+1},
\end{equation*}
we put ${\bf GL}_{\alpha}:=
\times_{i=1}^{d+1}{\bf GL}_{a_i}$ (we set ${\bf
GL}_0=\{e\}$). Let
\begin{equation*}{\rm Rep}(\mathcal
V_d, \alpha):={\rm Mat}_{a_1\times
a_2}\times\ldots\times{\rm Mat}_{a_1\times
a_{d+1}}\end{equation*} be the space of
$\alpha$-dimensional representations of ${\mathcal
V}_d$ endowed with the natural ${\bf
GL}_{\alpha}$-action   (we refer to \cite{kac1},
\cite{kac2}, \cite{kac3}, \cite{derksen-weyman}
\cite{sch1}, \cite{sch2} regarding the definitions and
notions of the representation theory of quivers). For
$\mathcal V_d$, the Euler inner product $\langle\ {,}\
\rangle$ on $\mathbb Z^{d+1}$ is given by
\begin{equation}\label{euler}
\langle(x_1,\ldots x_{d+1}), (y_1,\ldots,
y_{d+1})\rangle=(x_1y_1+\ldots
+x_{d+1}y_{d+1})-y_1(x_2+\ldots +x_{d+1}).
\end{equation}

It easily follows from the basic definitions, that the
following properties are equivalent:
 \begin{enumerate}
 \item[(i)] there is an open $G$-orbit in
 $(G/P_i)^d$;
 \item[(ii)] for $\gamma:=(l+1, i,\ldots, i)$,
 there is an open
${\bf GL}_\gamma$-orbit in ${\rm Rep}(\mathcal V_d,
\gamma)$.
\end{enumerate}

Note that by \cite[Corollary \,1 of
Proposition 4]{kac2}, (ii) is equivalent to the
following property
\begin{enumerate}
\item[(iii)] all the roots $\beta_i$ appearing in the
canonical decomposition of $\gamma$,
\begin{equation}\label{canon}
\gamma=\beta_1+\ldots+\beta_s,\end{equation} \noindent
are real, i.e., $\langle\beta_i, \beta_i\rangle=1$.
\end{enumerate}
Since there are combinatorial algorithms for finding
decomposition \eqref{canon}, see \cite{sch1},
\cite{sch2}, \cite{derksen-weyman} (the algorithm in
\cite{derksen-weyman} is very quick), (iii) and
 \eqref{euler} permit, in principle, to check
for every concrete $d, l, i$ whether  (i) holds or
not, and thereby to calculate ${\rm gtd}(G:G/P_i)$.
However we wish to obtain a closed formula for ${\rm
gtd}(G:G/P_i)$. We preface the corresponding argument
with the following observations.

\vskip 2mm

 (A) Let ${\mathcal V}_d^*$ be the
quiver  obtained from $\mathcal V_d$ by reversing the
directions of all the arrows. Let
\begin{equation*}
{\rm Rep}(\mathcal V_d^*, \alpha):={\rm
Mat}_{a_2\times a_1}\times\ldots\times{\rm
Mat}_{a_{d+1}\times a_{1}}
\end{equation*}
be the space of $\alpha$-dimensional representations
of ${\mathcal V}_d^*$ endowed with the natural ${\bf
GL}_{\alpha}$-action. The definition of ${\bf
GL}_{\alpha}$-actions readily shows that restricting
the map
\begin{equation*}
{\rm Rep}(\mathcal V_d, \alpha)\longrightarrow {\rm
Rep}(\mathcal V_d^*, \alpha),\quad
(A_1,\ldots,A_{d})\mapsto (A_1^{\top},\ldots,
A_d^{\top}),
\end{equation*}
to  any ${\bf GL}_\alpha$-orbit in ${\rm Rep}(\mathcal
V_d, \alpha)$ yields an isomorphism with some ${\bf
GL}_\alpha$-orbit in ${\rm Rep}(\mathcal V_d^*,
\alpha)$. Hence the existence of an open ${\bf
GL}_\alpha$-orbit in ${\rm Rep}(\mathcal V_d, \alpha)$
is equivalent to its existence in ${\rm Rep}(\mathcal
V_d^*, \alpha)$.

\vskip 2mm

(B) If either $a_1\leqslant a_i$ for every $i$, or
$a_1\geqslant a_2+\ldots + a_{d+1}$, then ${\rm
Rep}(\mathcal V_d, \alpha)$ contains an open ${\bf
GL}_\alpha$-orbit. This readily follows from the
definition of ${\bf GL}_\alpha$-action on ${\rm
Rep}(\mathcal V_d, \alpha)$.

\vskip 2mm

(C) Let $r_i$ be the $i$th fundamental reflection of
$\mathbb Z^{d+1}$, i.e.,
\begin{equation}\label{reflection}
r_i(\nu)= \nu-(\langle \nu, \alpha_i\rangle + \langle
\alpha_i, \nu\rangle)\alpha_i, \quad \nu\in \mathbb
Z^{d+1},\ \alpha_i=(0,\ldots,0,\underset{ \hat
i}1,0\ldots,0).
\end{equation}
It follows from \eqref{reflection},
\eqref{euler} that
\begin{gather}\label{formulae}
\begin{gathered}
r_i(\alpha)=\begin{cases}
(-a_1+a_2+\ldots+a_{d+1}, a_2,
\ldots, a_{d+1})&\text{for $i=1$},\\
(a_1,\ldots,a_{i-1}, a_1-a_i, a_{i+1},\ldots,
a_{d+1})& \text{for $i>1$.}
\end{cases}
\end{gathered}
\end{gather}
From \cite[2.3]{kac1}, \cite{kac3}, \cite[\S2,
Proposition 7]{sk}, and \eqref{formulae} we then
deduce the following.
\begin{enumerate}
\item[(C${}_1$)] Let $a_1\leqslant
a_2+\ldots+a_{d+1}$. Then ${\rm Rep}(\mathcal V_d,
\alpha)$ contains an open ${\bf GL}_\alpha$-orbit if
and only if ${\rm Rep}(\mathcal V_d^*, r_1(\alpha))$
contains an open ${\bf GL}_{r_1(\alpha)}$-orbit.
\item[(C${}_2$)] Let $a_1\geqslant a_i$ for all $i>1$.
Then ${\rm Rep}(\mathcal V_d, \alpha)$ contains an
open ${\bf GL}_\alpha$-orbit if and only if ${\rm
Rep}(\mathcal V_d^*, r_{d+1}\ldots r_2(\alpha))$
contains an open ${\bf GL}_{r_{d+1}\ldots
r_2(\alpha)}$-orbit.
\end{enumerate}

\vskip 1mm

We can now complete the proof of Theorem~{\rm
\ref{main3} for $G$ of type ${\sf A}_l$ using the
argument due to {\sc A.~Schofield}, \cite{sch3}.
Namely, we shall show that for every $\lambda\in
\mathfrak G:=
  \{(a_1,\ldots, a_{d+1})\in
  \mathbb N^{d+1}\mid a_2=\ldots =a_{d+1}\}$,
  we have
\begin{gather}\label{real}
\begin{gathered}
\text{${\rm Rep}({\mathcal V}_d, \lambda)$ contains an
open ${\bf GL}_{\lambda}$-orbit}\Longleftrightarrow
\langle \lambda, \lambda\rangle >0.
\end{gathered}
\end{gather}
By virtue of \eqref{Mm}, Definition~\ref{defgtdeg},
and \eqref{euler}, this claim immediately yields the
statement of Theorem~{\rm \ref{main3} for $G$ of type
${\sf A}_l$
 since for $\lambda=(n, a,\ldots, a)$, we have
\begin{equation}\label{length}
\langle \lambda, \lambda\rangle=n^2+da^2-nda.
\end{equation}

Turning to the proof of  claim, we call $\lambda_1$
and $\lambda_2\in \mathfrak G$ {\it congruent}
  if
 (i) $\langle \lambda_1,
\lambda_1\rangle= \langle \lambda_2,
\lambda_2\rangle$; (ii)
 ${\rm Rep}({\mathcal
V}_d, \lambda_1)$ contains an open ${\bf
GL}_{\lambda_1}$-orbit $\Longleftrightarrow$ ${\rm
Rep}({\mathcal V}_d, \lambda_2)$ contains an open
${\bf GL}_{\lambda_2}$-orbit. So, proving \eqref{real}
is equivalent to that with $\lambda$ replaced by
  a congruent vector. Recall that the quadratic
  form $\alpha\mapsto \langle\alpha, \alpha\rangle$
  on $\mathbb Z^{d+1}$ is invariant with
  respect to the group generated by
  $r_1,\ldots, r_{d+1}$.

Take $\lambda=(n, a,\ldots, a)\in \mathfrak G$. If
$a\geqslant n$, then by \eqref{length},
$\langle\lambda, \lambda\rangle=n^2+da(a-n)
>0$, and by (B),  ${\rm
Rep}({\mathcal V}_d, \lambda)$ contains an open ${\bf
GL}_\lambda$-orbit. This agrees with \eqref{real}.

Assume now that $a<n$. Then by (C${}_2$),
\eqref{formulae}, (A), vectors $\lambda$ and $(n, n-a,
\ldots, n-a)$ are congruent. Hence, in order to prove
\eqref{real} we may (and shall) assume that
\begin{equation}\label{ineq1}
n \geqslant 2a.\end{equation}

Consider now separately two cases: $2n\leqslant da$
and $2n> da$. If
\begin{equation}\label{ineq2}
2n\leqslant da,
\end{equation}
then from \eqref{euler}, \eqref{ineq1}, \eqref{ineq2}
we deduce
\begin{gather}\label{ineq3}
\begin{gathered}
\langle\lambda, \alpha_1\rangle + \langle \alpha_1, \lambda\rangle=2n-da\leqslant 0,\\
\langle\lambda, \alpha_i\rangle + \langle \alpha_i,
\lambda\rangle= 2a-n\leqslant 0,\quad i>1.
\end{gathered}
\end{gather}
The inequalities \eqref{ineq3} mean that $\lambda$
lies in the fundamental set $M$ of imaginary roots. By
\cite[Lemma 2.5]{kac1}, \cite[Corollary 1 of
Proposition 4]{kac2}, this yields that   there is no
open ${\bf GL}_\lambda$-orbit in ${\rm Rep}({\mathcal
V}_d, \lambda)$. On the other hand, $\lambda\in M$
yields $\langle\lambda, \lambda\rangle\leqslant 0$.
This agrees with~\eqref{real}.

Assume now that the second case holds, i.e.,
equivalently,
\begin{equation}\label{ineq4}
n> da-n.
\end{equation}
If $n\geqslant da$, then, by \eqref{length},
$\langle\lambda, \lambda\rangle=n(n-da)+da^2>0$, and,
by (B), ${\rm Rep}({\mathcal V}_d, \lambda)$ contains
an open ${\bf GL}_\lambda$-orbit. This agrees with
\eqref{real}. Assume now that $da>n$. Then, by
(C${}_1$), (A), \eqref{ineq4}, vectors $\lambda$ and
$(da-n, a,\ldots, a)$ are congruent. Thus, in view of
\eqref{ineq4}, proving \eqref{real} is reduced to that
with $\lambda=(n, a,\ldots, a)$ replaced by
$\lambda'=(n', a,\ldots, a)$ where $0<n'<n$. We then
can repeat the above arguments, starting from ``Take
$\lambda\ldots$'', with
 $\lambda$ replaced by $\lambda'$.
 Since the first coordinate of dimension
 vector strictly decreases
 via this process,
 the latter will eventually terminate. This completes the proof. {\hskip
8mm}~$\square$

\medskip


\noindent{\it Proof of Theorem~{\rm \ref{main1}.}}
Statements (i) and (ii) follow respectively from
statements (ii) and (iii) of Proposition~\ref{>1}.
Statement (iii) follows from Proposition~\ref{isogen},
and (iv)   from the Corollary of
Proposition~\ref{product}.

It is not difficult to deduce from \eqref{Mm} that
$m_{li}\leqslant l+2$ for every $i\leqslant l$, and
$m_{l1} =l+2$ (note that if $G={\bf SL}_{l+1}$, then
$G/P_1$ is ${\bf P}^{l}$ endowed with the natural
${\bf SL}_{l+1}$-action, and by Theorem~\ref{main3},
equality $m_{l1} =l+2$ expresses the well known
elementary fact that this action is generically
$(l+2)$-transitive). Statement (v) now follows from
Definition~\ref{defgtdeg}, Theorem~\ref{main2}, and
Theorem~\ref{main3}. {\hskip 8mm}~$\square$

\providecommand{\bysame}{\leavevmode\hbox
to3em{\hrulefill}\thinspace}


\begin{thebibliography}{XXXXx}

\bibitem[ABS]{abs}
  {\sc  H.~Azad, M.~Barry, G.~Seitz},
  \emph{On the structure of parabolic
  subgroups}, Comm. in Algebra
  \textbf{18(2)} (1990), 551--562.

  \bibitem[Bor]{borel}
  {\sc  A.~Borel},
  {\it Linear Algebraic
  Groups}:
  Second Enlarged Edition,
  Graduate Texts in
  Mathematics,
  Vol. 126, Springer-Verlag, 1991.

  \bibitem[Bou]{bourbaki}
  {\sc  N.~Bourbaki},
  \emph{Groupes et alg\`ebres de Lie},
  Chap. IV, V,  VI, Hermann, Paris, 1968.

  \bibitem[DW]{derksen-weyman}
  {\sc H.~Derksen, J.~Weyman},
  \emph{On the canonical decomposition of quiver
  representations},
  Compositio Math. {\bf 133} (2002), 245--265.

  \bibitem[Hu]{humphreys}
  {\sc J.~E.~Humphreys},
  \emph{Linear Algebraic Groups},
  Springer-Verlag, New York, Heidelberg,
  Berlin, 1975.

  \bibitem[Ka${}_1$]{kac1}
  {\sc V.~Kac},
  \emph{Infinite root systems,
  representations of graphs and invariant
  theory}, Invent. Math. {\bf 56}
  (1980), 57--92.

  \bibitem[Ka${}_2$]{kac2}
  {\sc V.~Kac},
  \emph{Infinite root systems,
  representations of graphs and invariant
  theory} II, J. Algebra {\bf 78}
  (1982), 141--162.

  \bibitem[Ka${}_3$]{kac3}
  {\sc V.~Kac},
  \emph{Root systems, representations of quivers, and
  Invariant theory}, Lect. Notes Math., Vol. 996,
  Springer-Verlag, 1983, 73--108.


\bibitem[Kime${}_1$]{kime1}
  {\sc  B.~Kimel'fel'd},
  \emph{Reductive groups which are locally transitive
  on flag manifolds of orthogonal groups}, Tr. Tbilis.
  Mat. Inst. Razmadze
  \textbf{LXII} (1979), 49--75
  (in Russian). English transl.:
  Sel. Math. Sov. {\bf 4} (1985),
  107--130.

  \bibitem[Kime${}_2$]{kime2}
  {\sc  B.~Kimel'fel'd},
  \emph{Homogeneous domains
  on flag manifolds}, J. Math. Anal. Appl.
  \textbf{121} (1987), 506--588.

  \bibitem[Kimu]
  {kimu1}
  {\sc  T.~Kimura},
  \emph{A classification of prehomogeneous
  vector spaces of simple algebraic
  groups with scalar multiplications},
 J. Algebra {\bf 83} (1983), 72--100.





\bibitem[KKIY]{kkiy}
  {\sc  T.~Kimura, S.-i.~Kasai,
  M.~Inuzuka, O.~Yasukura},
  \emph{A classification of} $2$-{\it simple
  prehomogeneous vector spaces of type} I,
  J.~Algebra
  \textbf{114} (1988), 369--400.

  \bibitem[KKTI]{kkti}
  {\sc  T.~Kimura, S.-i.~Kasai,
  M.~Taguchi, M.~Inuzuka},
  \emph{Some P.V.-equivalences and
  a classification of} $2$-{\it simple
  prehomogeneous vector spaces of type} II,
  Trans.~AMS
  \textbf{308} (1988), No. 2,
  433--494.

\bibitem[Kn]{knop}
  {\sc  F.~Knop},
  \emph{Mehrfach transitive Operationen algebraischer
  Gruppen}, Arch. Math.
  \textbf{41} (1983), 438--446.

  \bibitem[Li]{littelmann}
{\sc P.~Littelmann}, \emph{On spherical double cones},
J. Algebra {\bf 166} (1994), no. 1, 142--157.

  \bibitem[Lu]{luna}
  {\sc  D.~Luna},
  \emph{Sur les orbites ferm\'ees des groupes
  alg\'ebriques reductifs}, Invent. Math.
  \textbf{16} (1972), 1--5.

  \bibitem[Ni]{nisn} {\russc E.~A.~Nisnevi{ch}},
  {\rusit Perese{ch}enie podgrupp reduktivnyh grupp i
  stabil{p1}nost{p1} de{i0}stvi{ya}},
  {\rus Dokl. Akad. Nauk BSSR}
   {\bf 17} (1973),
  785--787.
  ({\sc E.~A.~Nisnevich},
  \emph{Intersection of subgroups of reductive
  groups and stability of action},
  Dokl. Akad. Nauk BSSR {\bf 17} (1973),
  785--787 (in Russian)).

  \bibitem[OV]{oni-vinb} {\russc E.~B.~Vinberg,
  A.~L.Oniwik}, {\rusit Seminar po gruppam Li i
  algebrai{ch}eskim gruppam}, {\rus Nauka, Moskva},
  1988. Engl. transl.:
  {\sc A.~L.~Onishchik, E.~B.~Vinberg},
  \emph{Lie Groups and Algebraic Groups},
  Springer-Verlag, Berlin, Heidelberg, New York,
  1990.

  \bibitem[P]{popov}
  {\sc V.~L.~Popov}, \emph{Tensor
  product decompositions and products
  of flag varieties}, preprint, Steklov Mathematical Institute, Russian Academy of Sciences, 2005.

  \bibitem[PV${}_1$]{popov-vinberg1}
  {\russc E.~B.~Vinberg,
  V.~L.~Popov}, {\rusit Ob odnom klasse
  kvaziodnorodnyh affinnyh mnogoobrazi{i0}},
  {\rus  Izv. AN SSSR, ser. mat.} {\bf 36} (1972),
  749--763. Engl. transl.:
  {\sc  E.~B.~Vinberg, V.~L.~Popov},
  \emph{On a class of quasihomogeneous affine
  varieties},
  Math. USSR, Izv. {\bf 6} (1973), 743--758.

  \bibitem[PV${}_2$]{popov-vinberg2}
  {\russc E.~B.~Vinberg,
  V.~L.~Popov}, {\rusit Teori{ya} invariantov},
   {\rus So\-vr. probl. matem.,
Fund. na\-prav\-l.}, {\rus VI\-NI\-TI, Moskva, t. 55},
1989, {\rus str.} 137--309. Engl. transl.:
  {\sc V.~L.~Popov, E.~B.~Vinberg},
  \emph{Invariant Theory}, Encycl. of Math. Sci., Vol.
  55, Springer-Verlag, Heidelberg, 1994, pp. 123--284.

\bibitem[Ri${}_1$]{richardson0}
{\sc R.~Richardson}, \emph{Conjugacy classes in
parabolic subgroups of semisimple algebraic groups},
Bull. London Math. Soc. {\bf 6} (1974), 21--24.


\bibitem[Ri${}_2$]{richardson}
{\sc R.~Richardson}, \emph{Finiteness theorems for
orbits of algebraic groups}, Indag. Math. {\bf 88}
(1985), 337--344.

\bibitem[Ro]{rosenlicht}
{\sc M.~Rosenlicht}, \emph{A remark on quotient
spaces}, An. Acad. Brasil. Ci. {\bf 35} (1963),
487--489.


\bibitem[RRS]{rrs}
{\sc R.~Richardson, G.~R\"ohrle, R.~Steinberg},
\emph{Parabolic subgroups with abelian unipotent
radical}, Invent. math. {\bf 110} (1992), 649--671.

\bibitem[R\"oh]{roh}
{\sc G.~R\"ohrle}, \emph{On the structure of parabolic
subgroups in algebraic groups}, J. Algebra {\bf 157}
(1993), 80--115.


  \bibitem[SK]{sk}
  {\sc  M.~Sato, T.~Kimura},
  \emph{A calssification of irreducible prehomogeneous
  vector spaces and their
  relative invariants}, Nagoya Math. J.
  \textbf{65} (1977), 1--155.


  \bibitem[Sch${}_1$]{sch1}
  {\sc  A.~Schofield}, \emph{General representations
  of quivers}, Proc. London Math. Soc. {\bf 65} (1992),
  46--64.

\bibitem[Sch${}_2$]{sch2}
  {\sc  A.~Schofield},
  \emph{Birational classification
  of moduli spaces of vector bundles
  over $\mathbb P^2$},
  Indag. Math. (N.S.) {\bf 12} (2001), no. 3,
  433--448.


  \bibitem[Sch${}_3$]{sch3}
  {\sc  A.~Schofield}, \emph{Letter to V.~L.~Popov},
  January 31, 2005.


\bibitem[Sp]{springer}
  {\sc T.~A.~Springer},
  \emph{Linear Algebraic Groups. Second
  Edition}, Progress in Math., Vol. 9,
  Birkh\"auser, Boston, Basel, Berlin,
  1998.

  \bibitem[VK]{vk} {\russc E.~B.~Vinberg,
  B.~N.~Kimel{p1}fel{p1}d}, {\rusit Odnorodnye oblasti
  na flagovyh mnogoobrazi{ya}h i sferi{ch}eskie
  podgruppy
  poluprostyh grupp Li}, {\rus Funkc. anal. i ego
  pril.} {\bf 12} (1978), 12--19. Engl. transl.:
  {\sc  E.~B.~Vinberg, B.~N.~Kimelfeld},
\emph{Homogeneous domains on flag manifolds
  and spherical subgroups of semisimple Lie groups}, Funct. Anal.
  Appl.
  \textbf{12} (1979), 168--174.



\end{thebibliography}
\end{document}